\documentclass{article}
\usepackage{fullpage}
\usepackage{hyperref}

\usepackage{amsmath,amssymb,amsfonts,mathrsfs,amsthm}
\usepackage{array}
\usepackage[utf8]{inputenc}
\usepackage{listings}
\usepackage{mathtools}
\usepackage{pdfpages}
\usepackage[textsize=footnotesize,color=green]{todonotes}
\usepackage{bm}
\usepackage{tikz}
\usepackage[normalem]{ulem}
\usepackage{hhline}

\usepackage{algorithm}
\usepackage[noend]{algpseudocode}
\usepackage{algorithmicx}

\usepackage{graphicx}
\usepackage{subfigure}

\usepackage{color}
\usepackage{pdflscape}
\usepackage{pifont}

\newcommand{\vect}[1]{\ensuremath\boldsymbol{#1}}

\newcommand{\pd}[2]{\frac{\partial#1}{\partial#2}}

\newcommand{\nor}[1]{\left\| #1 \right\|}

\newcommand{\LRp}[1]{\left( #1 \right)}

\newcommand{\LRb}[1]{\left| #1 \right|}

\newcommand{\eval}[2][\right]{\relax
  \ifx#1\right\relax \left.\fi#2#1\rvert}

\def\vectwo#1#2{\left[
\begin{array}{c}
#1\\
#2\\
\end{array}
\right]}

\title{A Comparison of High-Order Interpolation Nodes for the Pyramid}

\author{Jesse Chan\thanks{Department of Computational and Applied Mathematics, Rice University, 6100 Main St, Houston, TX, 77005} \and T. Warburton\footnotemark[1]}

\date{}
\begin{document}

\maketitle

\begin{abstract}
The use of pyramid elements is crucial to the construction of efficient hex-dominant meshes \cite{bergot2010higher}.  For conforming nodal finite element methods with mixed element types, it is advantageous for nodal distributions on the faces of the pyramid to match those on the faces and edges of hexahedra and tetrahedra.  We adapt existing procedures for constructing optimized tetrahedral nodal sets for high order interpolation to the pyramid with constrained face nodes, including two generalizations of the explicit Warp and Blend construction of nodes on the tetrahedron \cite{warburton2006explicit}.  

\end{abstract}



\section{Introduction}

In recent years, high order finite element methods have been shown to have significant advantages over low order finite elements in a variety of areas.  In particular, for smooth solutions, they tend to converge more rapidly under both order and mesh refinement, and for wave propagation, display less numerical dissipation under time-marching schemes than low order discretizations \cite{deville2002high}.  High order finite element stifness and mass matrices are also typically block structured, allowing for efficient local computations in matrix-free applications.  

Quadrilateral and hexahedral elements may offer significant benefits over triangular and tetrahedral elements in high order finite element methods as well --- exploiting a tensor-product structure allows for simplified data structures, as well as fast, low-memory applications of operators in matrix-free methods.  Additionally, high order discontinuous Galerkin (DG) methods benefit from the use of hexahedral elements by reducing the number of flux calculations and also reducing the number of degrees of freedom compared to tetrahedral elements, while still delivering the same global order of approximation.  However, while tetrahedral mesh generation has developed to be able to mesh near-arbitrary geometries, it is difficult to mesh arbitrary geometries using purely hexahedral meshes.  Often, it is at most possible to construct hex-dominant meshes, which contain primarily hexahedral elements, but also tetrahedral, wedge, and pyramid elements  \cite{Demkowicz:2007:CHF:1564840, bergot2010higher, baudouin2014frontal}.

We are interested in developing $H^1$-conforming finite element methods using nodal (Lagrange) basis functions.  For ease of conformity between elements, it is preferable to have the same distribution of nodes on each type of element face.  The nodal distribution may also be chosen in such a way that it minimizes the interpolation error for the given nodal basis via minimization of the Lebesgue constant.  Defining such distributions for quadrilateral and hexahedral elements is typically done through a tensor-product construction, while construction of optimal distributions for triangular and tetrahedral elements is well-explored \cite{roth2005nodal, rapetti2012generation, chen1996optimal, hesthaven1998electrostatics, warburton2006explicit, taylor2000algorithm}.  For the wedge element, a simple tensor-product extrusion of nodes on a triangular face allows for conformity with both tetrahedral and hexahedral elements with a reasonable Lebesgue constant.   However, for the pyramid, optimized distributions of nodes have received less attention \cite{de2013polynomial, bos2011geometric, bergot2010higher, gassner2009polymorphic}, and are the focus of this article.  


\section{Finite element spaces and orthogonal bases on pyramids}
\label{sec:bases}
Since pyramid and wedge elements are used primarily to link tetrahedral, hexahedral, and wedge elements together, it is appropriate to require that the trace space on triangular faces of the pyramid and tetrahedra are the same, and similarly for the quadrilateral faces of the pyramid and hexahedra.  Additionally, we require that the edge trace spaces of the pyramid, tetrahedra, and hexahdrea are the same.    

Standard trace spaces for tetrahedral and hexahedral finite element spaces contain polynomials of total degree or independent degree $N$ respectively; however, it has been shown by Bedrossian \cite{bedrosian1992shape}, Wieners \cite{wieners1997conforming} and Nigam and Phillips \cite{nigam2012high} that it is impossible to construct a basis with polynomial trace spaces on a pyramid with purely polynomial functions.  Consequentially, conforming finite element spaces for pyramids necessarily contain rational functions.

Nigam and Phillips \cite{nigam2012high} give a comprehensive construction of arbitrary order $H^1$, $H({\rm curl})$, $H({\rm div})$ and $L^2$ finite element spaces for pyramids in terms of a monomial basis.  
Prior to this, $L^2$ orthogonal basis functions for the pyramid were used early on in $hp$ and spectral finite element simulations by Warburton, Sherwin, and Karniadakis in \cite{warburton1999spectral, sherwin1998spectral, karniadakis2013spectral}, and partially orthogonalized high order finite element spaces were used in $hp$-adaptive conforming finite elements for the Maxwell's equations with exact sequence by Zaglmayr \cite{zaglmayr} and Demkowicz et al \cite{Demkowicz:2007:CHF:1564840}.  

Bergot, Cohen, and Durufle present an alternative $L^2$ orthogonal basis of optimal dimension on the pyramid in \cite{bergot2010higher}.  Comparisons are made with existing finite element spaces, and the basis of Warburton, Sherwin, and Karniadakis is shown to be suboptimal for higher order, while the basis of Nigam and Phillips is shown to be optimal, but of larger dimension than necessary\footnote{A new basis with smaller dimension is presented by Nigam and Phillips in \cite{nigam2012numerical}. The resulting approximation space matches that of Bergot, Cohen, and Durufle.  It should be noted that Zaglmayr's basis spans the same approximation space as well.}.  Though the basis itself is rational, its traces are polynomial, satisfying our requirement that the trace spaces match with tetrahedral and hexahedral elements.  We will use the basis of Bergot, Cohen, and Durufle in our construction of nodal basis functions.  For a reference pyramid (see Figure~\ref{fig:split}) with coordinates $r, s, t$ such that
\[
r,s \in [t-1,1-t]  \quad t \in [0,1],
\]
we define the basis functions $P_{ijk}(r,s,t)$ as
\[
P_{ijk}(r,s,t) = P_i^{0,0}\LRp{\frac{r}{1-t}}P_j^{0,0}\LRp{\frac{s}{1-t}}(1-t)^{c}P_k^{2(c+1),0}\LRp{2t-1},
\]
where 
\[
c = \max(i,j), \quad 0\leq i,j \leq N, \quad 0 \leq k \leq N-c.
\]
and $P^{a,b}_n(x)$ is the Jacobi polynomial of order $n$, orthogonal with respect to the weight $(1-x)^a(1+x)^b$.  For a given $N$, the above procedure produces $N_p$ distinct basis functions, where 
\[
N_p = \frac{(N+1)(N+2)(2N+3)}{6}
\]
We may order these orthogonal functions arbitrarily as $\phi_j(r,s,t)$ from $j = 1, \ldots, N_p$.  

Using the above basis functions, we may define the generalized Vandermonde matrix 
\[
V_{ij} = \phi_j(r_i,s_i,t_i), \quad i = 1, 2, \ldots, N_p
\]
where $(r_i, s_i, t_i)$ are a set of nodal points contained inside or on the boundary of the reference pyramid.  Then, a nodal basis may be constructed using elements of the inverse of the Vandermonde matrix 
\[
\ell_i(r,s,t) = \sum_{j = 1}^{N_p}(V^{-1})_{ij}\phi_j(r,s,t).
\]


\section{Construction of optimized nodal sets}

In the discrete setting, enforcing conformity for nodal finite element methods reduces to the matching of nodal degrees of freedom from one element face to another element face.  This motivates the requirement that the distribution of surface nodes for both pyramids and wedges should match the distribution of face nodes on tetrahedra and hexahedra.  Furthermore, it is assumed that the nodes on the edges of each element follow a Gauss Legendre-Lobatto (GLL) distribution.  

For quadrilateral and hexahedral elements, the simplest construction is to take a tensor product of 1D GLL nodes.  This choice of nodes is central in the formulation of Spectral Element Methods (SEM); combined with GLL quadrature and tensor product evaluation of operators, the result is a highly efficient method for tensor product meshes \cite{patera1984spectral, deville2002high, karniadakis2013spectral}.  For wedges, the nodal distribution may also be easily constructed --- assuming a given distribution of face nodes for the tetrahedra with GLL nodes on the edges, we can extrude this distribution in the direction orthogonal to the face, which results in a tensor product GLL structure on the quadrilateral face of the wedge.  

For pyramids, the ideal nodal distribution is less clear.  It is possible to simply combine an arbitrary surface and interior node distribution; for example, both Bergot, Cohen, and Durufle and Gassner et al combine the electrostatic nodes of Hesthaven on the triangular faces of the pyramid with an appropriate number of tensor product GLL nodes in the interior and on the base, similar to the construction of Stroud quadrature \cite{bergot2010higher, gassner2009polymorphic}.  However, it is possible to choose a more tailored distribution of interior nodes in order to optimize some measure of quality of the nodal set.  

\subsection{Metrics of quality and optimization strategies}

A common aim in the construction of nodal sets is to minimize the Lebesgue constant the nodal distribution.  Given a set of $N_p$ nodal points $\{x_1,\ldots, x_{N_p}\}$, the Lebesgue constant is defined in terms of the Lagrange interpolatory basis functions $\ell_i(\vect{x})$ as
\[
\Lambda = \max_{\vect{x} \in K}\sum_{i = 1}^{N_p} \LRb{\ell_i(\vect{x})}.
\]
The Lebesgue constant bounds the interpolation error in the max norm, such that
\[
\|f-f_N\|_{\infty} \leq \LRp{\Lambda + 1} \|f - f_N^*\|_{\infty},
\]
where $f_N$ is the order $N$ interpolant of $f$, $f_N^*$ is the best order $N$ approximation to $f$, and $\Lambda$ is the Lebesgue constant for the given interpolation points.

For the tetrahedron, nodal distributions which minimize the Lebesgue constant have been explored in great detail.  Chen and Bab\u{u}ska minimized the $L^2$ norm of the sum of the Lagrange basis --- an $L^2$ analogue of the Lebesgue constant --- in \cite{chen1996optimal}.  Hesthaven constructed points using an analogy to electrostatics; a nodal distribution is determined by finding stationary distributions of charges, which are in turn related to the zeros of specific Jacobi polynomials \cite{hesthaven1998electrostatics}.  

The Warp and Blend construction of nodes on triangles and tetrahedra involve defining a warp, or a displacement which shifts equispaced nodes to GLL nodes on an edge, and blending this warping function into the interior of faces and volumes using techniques from curvilinear mesh generation \cite{warburton2006explicit, hesthaven2007nodal}.  The subsequent blends are linear in each barycentric coordinate, but an additional scaled quadratic warp term can be introduced in order to further optimize the Lebesgue constant over the resulting distribution of points.  

Another metric for quality of a nodal distribution is the determinant of the Vandermonde matrix, and distributions which maximize the determinant of the Vandermonde matrix are referred to as Fekete nodes.  
The maximization of the determinant is also an attractive alternative to direct optimization of the Lebesgue constant since there exist analytic expressions for the determinant of the Vandermonde matrix, which can be exploited in algorithms to compute Fekete nodes.
  For example, Taylor, Wingate and Vincent utilized analytic expressions for the determinant to develop a pseudo-time steepest ascent algorithm which can be used to iteratively move from an arbitrary node set to the Fekete node set of a given order \cite{taylor2000algorithm}.  

Other work has explored the optimization of interpolation points for more general shapes; Gassner et al \cite{gassner2009polymorphic} give a general construction for a polygon using a barycentric mapping, while \cite{2014arXiv1407.3291G, van2014approximating, narayan2013constructing} give algorithms for the approximation of good interpolation points on general domains in multiple dimensions.  Bos, Calvi, Levenberg, and Vianello construct nodal points using the concept of Weakly Admissible Meshes (WAM) \cite{bos2011geometric}, which are sequences of subsets $A_N$ of an element over which an upper bounds
\[
\|p(x)\|_{\infty} \leq C(A_N) \|p(x)\|_{\infty, A_N}, \quad \forall p\in P_N^d
\]
can be shown.  For Fekete points generated from weakly admissible meshes, $C(A_N)$ can be used to bound the Lebesgue constant, though the above bound only holds for polynomials of order $N$, while conforming finite element spaces on pyramids contain rational functions.  

A more recent innovation is the application of numerical linear algebra techniques to compute sets of so-called ``approximate Fekete'' points \cite{de2013polynomial}.  Bos et al \cite{bos2010computing} and Sommariva et al \cite{ Sommariva20091324} utilized similar ideas in the context of numerical linear algebra to compute ``approximate Fekete'' points.  We note that these points do not refer to the approximations of Fekete points detailed by Taylor, Wingate, and Vincent, and instead relate to a discrete equivalent of Fekete points, which use the concept of WAM to characterize ``approximate Fekete'' and ``discrete Leja'' points.  These points may then be computed using manipulations and factorizations of a Vandermonde matrix.  

\subsection{Face nodal distributions}

For the following experiments, we enforce a fixed distribution of nodes on the faces of the pyramid.  For quadrilateral faces, we choose SEM/tensor product GLL nodes for conformity with hexahedral elements.  
For triangular faces of the pyramid, we fix the face nodes to be identical to those of the Warp and Blend tetrahedra.  As noted above, the choice of nodal distributions for the tetrahedron is less obvious.  For consistency in the comparison of different pyramid nodal sets, we fix the triangular faces to be the Warp and Blend nodes for the tetrahedra, which give competitive Lebesgue constants, but maintain a simple and explicit construction for any order $N$.  

We note that it may be possible to optimize nodal distributions over triangular faces to minimize the Lebesgue constants for tetrahedra, wedge, and pyramid elements simultaneously.  
The work described here does not consider this option, though we hope to explore this in the future. 

\subsection{Fekete points}

Fekete points may be approximately constructed on the pyramid using the steepest ascent procedure from Taylor, Wingate, and Vincent, where nodal positions $r_i, s_i, t_i$ are taken to be stationary distributions of the ODEs
\[
\pd{r_i}{t}{} = \pd{\ell_i}{r}{},\qquad \pd{s_i}{t}{} = \pd{\ell_i}{s}{}, \qquad \pd{t_i}{t}{} = \pd{\ell_i}{t}{}.
\]
\begin{algorithm}
\begin{algorithmic}[1]
\Procedure{ODE steepest ascent Fekete}{}
\State Initialize positions of points $r_i, s_i, t_i$.
\State Initialize $t = 0$, timestep $dt$, and tolerance $tol$.  
\While{$t > 0$ and $\max_i \LRb{\pd{r_i}{t}{}} > tol$, $\max_i \LRb{\pd{s_i}{t}{}} > tol$, $\max_i \LRb{\pd{t_i}{t}{}} > tol$}
\State Compute, at current positions $r_i, s_i, t_i$,
\[
\pd{r_i}{t}{} = \pd{\ell_i}{r}{}, \quad \pd{s_i}{t}{} = \pd{\ell_i}{s}{}, \quad \pd{t_i}{t}{} = \pd{\ell_i}{t}{}
\]   
\State Time-march using LSERK-4 to determine new positions $r_i, s_i, t_i$.
\EndWhile
\State \textbf{return:} $r_i, s_i, t_i$.
\EndProcedure
\end{algorithmic}
\caption{ODE steepest ascent method for approximation of Fekete nodes \cite{taylor2000algorithm}.}
\label{alg:ode}
\end{algorithm}

The steepest ascent procedure in Algorithm~\ref{alg:ode} is known to be sensitive to initial conditions \cite{taylor2000algorithm}, so we initialize the nodal positions at time $t=0$ to a reasonable initial distribution by using the Warp and Blend points defined in Section~\ref{sec:pyramidWB}.  We solve the resulting system of ODEs using fourth order Runge-Kutta, and terminate the steepest ascent procedure when the maximum change over all nodes in their $r$, $s$, or $t$ positions is less than ${10^{-10}}/{dt}$ over a timestep.  Additionally, we do not update surface node positions, which fixes face distributions for conformity with other elements, which is done by setting the velocities $\pd{r_i}{t}{}$, $\pd{s_i}{t}{}$, and $\pd{t_i}{t}{} = 0$ for all nodes on the the pyramid faces.  The resulting nodes are referred to as ``Fekete'' in the numerical results.  

\begin{figure}
\centering
\subfigure{\includegraphics[width=.475\textwidth]{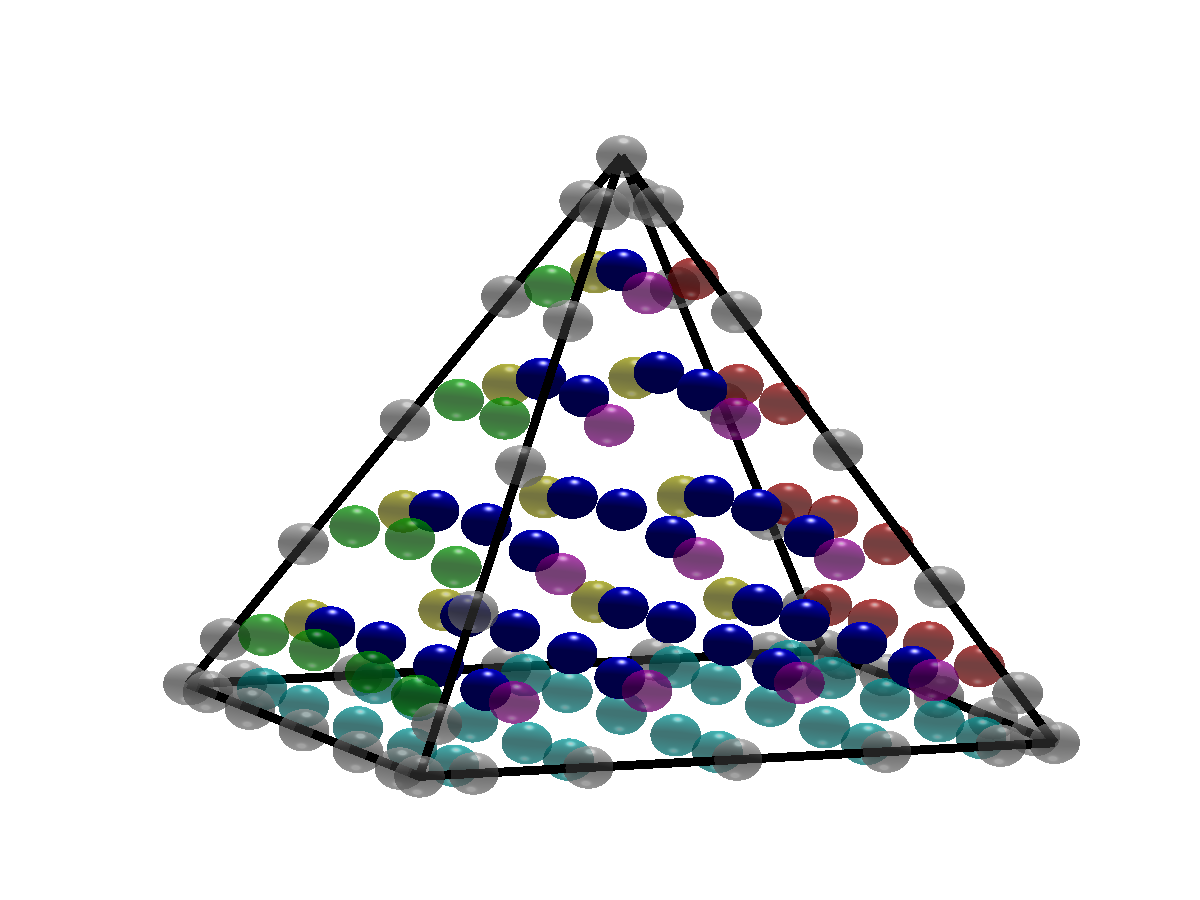}}
\subfigure{\includegraphics[width=.475\textwidth]{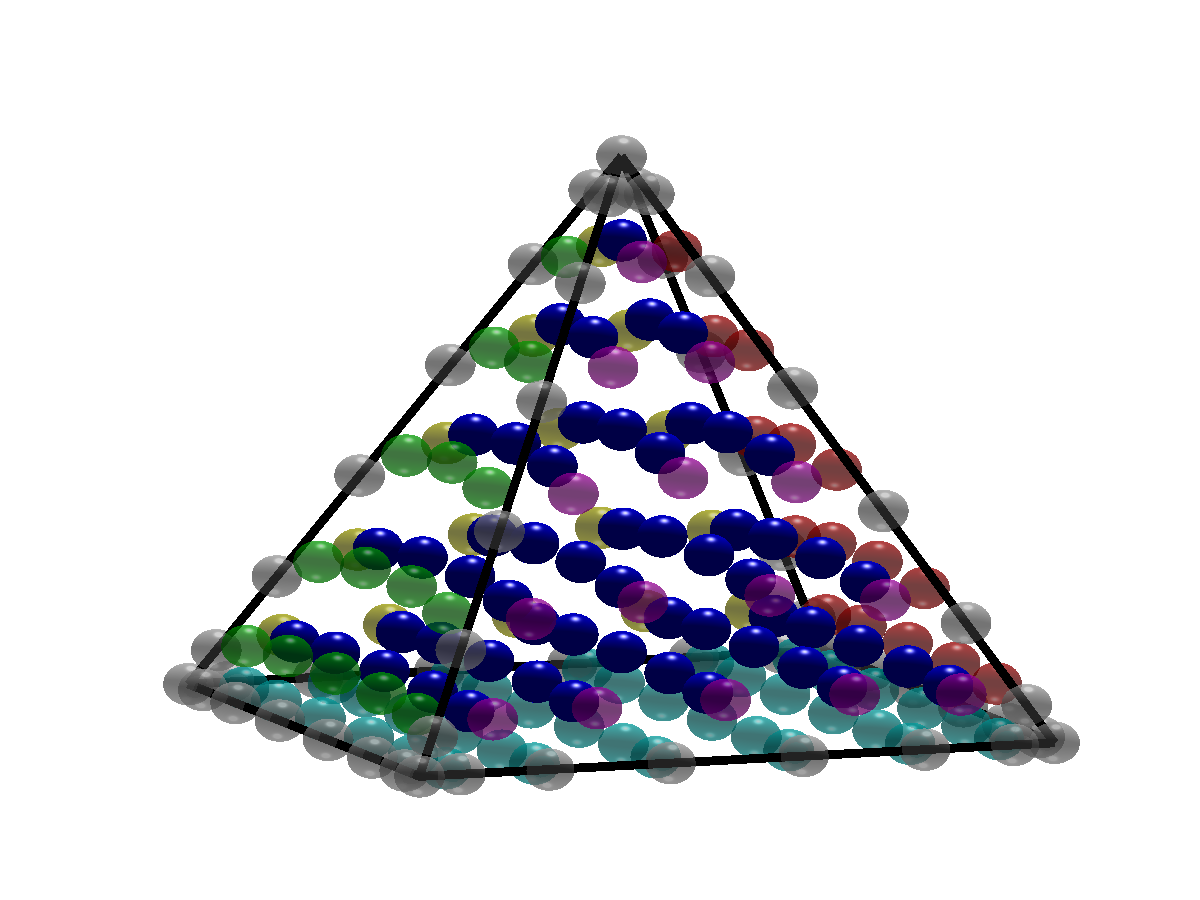}}
\caption{Fekete nodes for $N = 6$ (left) and $N=7$ (right) nodes for the pyramid.  Surface nodes (which have fixed nodal distributions) are shown as transparent spheres.  }
\label{fig:pyrFek}
\end{figure}

\subsection{``Approximate Fekete'' points}

Procedures to compute ``approximate Fekete'' points are given in \cite{bos2010computing, Sommariva20091324}.  An advantage of these algorithms apart from their relatively fast speed is their flexibility; since the algorithms are built around tools from numerical linear algebra, they are generalizable to a larger variety of domains and bases.  Most of these algorithms require only evaluations of some linearly independent basis at arbitrary points with which to generate a Vandermonde matrix, and are directly applicable to the pyramid.  We consider two procedures, both of which may be used to determine a distribution of interpolation points lying inside the pyramid.  

The first of these procedures is a greedy algorithm for maximizing the determinant of the Vandermonde matrix.  Supposing that $K$ is the pyramid, given a set of $N_s$ sample points $\{\vect{x}_i\}_{i = 1}^{N_s} \in K$, the following greedy algorithm chooses $N_p$ indices $i_k$ such that $\vect{x}_{i_1},\ldots, \vect{x}_{i_{N_p}}$ maximize the determinant of the Vandermonde matrix.  The algorithm selects these points sequentially based on the column with maximum norm.  The matrix is then updated by projecting out the component of the maximum norm column from all other columns.  The procedure is detailed in Algorithm~\ref{alg:greedy}.

\begin{algorithm}
\begin{algorithmic}[1]
\Procedure{Greedy ``approximate Fekete''}{}
\State Select points $\vect{x}_1,\ldots, \vect{x}_{N_s} \in K$.  
\State Construct $V_{ij} = \phi_j(\vect{x}_i)$, normalize the columns.  
\For{$k = 1,\ldots, N_p$}
\State Choose $i_k = {\rm arg max}_i \nor{V(:,i)}$.  
\For{$j \neq i_k$}
\State Orthogonalize $V(:,j)$ with respect to $V(:,i_k)$.  
\EndFor
\EndFor
\State \textbf{return:} $\vect{x}_{i_1}, \ldots, \vect{x}_{i_{N_p}}$.
\EndProcedure
\end{algorithmic}
\caption{Greedy selection of ``approximate Fekete'' points from a sample set \cite{Sommariva20091324}.}
\label{alg:greedy}
\end{algorithm}

Since we are interested in enforcing conformity between elements, we may introduce small modifications to fix the nodal distribution on pyramid faces.  Suppose that we wish to force inclusion of the first $N_b$ points $\vect{x}_1,\ldots, \vect{x}_{N_b}$ in the index set $\{i_k\}_{k = 1}^{N_p}$; we may then skip step 1 of the above algorithm for $k \leq N_b$. For the pyramid, we may force the inclusion of predetermined face nodes into the ``approximate Fekete'' set using this modification.   
We implement the above algorithm to determine ``approximate Fekete'' nodes on the pyramid with constrained face distributions, and refer to the resulting points as ``Greedy'' in the numerical results.  

The second procedure we consider is an iterative refinement method to determine ``approximate Fekete'' nodes based on the QR decomposition proposed by Sommariva and Vianello.  Suppose that $N_s$ denotes the number of sample points and $N_p$ the dimension of an arbitrary basis $\{\phi_1,\ldots, \phi_{N_p}\}$, initialize $V_0$ as the $N_s \times N_p$ Vandermonde matrix and $P_0$ as the $N_p\times N_p$ identity matrix.  Algorithm~\ref{alg:QR} details this process:
\begin{algorithm}
\begin{algorithmic}[1]
\Procedure{Iterative refinement ``approximate Fekete''}{}
\State Select points $\vect{x}_1,\ldots, \vect{x}_{N_s} \in K$, construct $V_{ij} = \phi_j(\vect{x}_i)$.  
\State Initialize $V_0 = V$, $P_0 = I$.
\For{$k = 0,\ldots, s-1$}
\State Compute the QR decomposition of $V_k = Q_k R_k$, and set
\[
V_{k+1} = V_kR_k^{-1}, \quad P_{k+1} = P_kR_k^{-1}
\]
\EndFor
\State Select $m_j \neq 0$.
\State Set $\mu = (P_s^T)^{-1}m, \quad w = (V_s^T)^{-1}\mu$.
\State Select $N_p$ points $i_k$ such that $w_{i_k} \neq 0$.  
\State \textbf{return:} $\vect{x}_{i_1}, \ldots, \vect{x}_{i_{N_p}}$.
\EndProcedure
\end{algorithmic}
\caption{Iterative refinement selection of ``approximate Fekete'' points from a sample set \cite{Sommariva20091324}.  }\label{alg:QR}
\end{algorithm}

The iterative refinement method has the advantage of being broadly applicable; the basis may be arbitrary, and the initial Vandermonde matrix $V_0$ may be low rank, since the construction of $V_s$ in the iterative refinement step will produce a nonsingular matrix.  

If we wish to enforce a fixed distribution of nodes on the faces of the pyramid, we may modify the basis in the above algorithm to be zero on over all nodes on the faces.  We do so by taking the underdetermined Vandermonde matrix $V_b$
\[
V_{b,ij} = \phi_j(x_i), \quad j = 1,\ldots, N_p, \quad i = 1,\ldots, N_b
\]
where $N_p$ is the dimension of the basis $N_b < N_p$ is the number of points on the surface/boundary of the pyramid. We may then take the $N_p-N_b$ linearly dependent combinations of columns of this matrix corresponding as coefficients which define functions that are zero on the faces --- in other words, an interior basis.\footnote{Constructing $V_i$ from a basis consisting purely of bubble functions supported only in the interior of the pyramid is another possibility.}  

Using $V_i$ in lieu of $V_0$ in the above algorithm allows us to compute ``approximate Fekete'' points only in the interior of the pyramid.  We refer to these resulting points as ``QR'' in the numerical results.  

For both ``approximate Fekete'' procedures, we take the sample set for a given $N$ to be Stroud-style equispaced points in the pyramid with $N^2+1$ points per edge.  This choice is motivated by a theorem of Bos and Levenberg which states that, in 1D, ``approximate Fekete'' points have the same asymptotic distribution as those of the true Fekete points for this choice of sample points \cite{bos2010computing}.  Further increasing the number of sampling points was not significantly correlated with an improvement in the Lebesgue constant of the resulting ``approximate Fekete'' point set.

\subsection{Pyramid nodes by Warp and Blend}

For conformity, we  assume that the nodal distribution on both triangular and quadrilateral faces is fixed.  For quadrilateral faces, we will choose SEM nodes for conformity with quadrilateral elements, and for triangular faces, we will choose the Warp and Blend optimized points of Warburton \cite{warburton2006explicit}.  

\begin{figure}
\centering
\subfigure{\includegraphics[width=.475\textwidth]{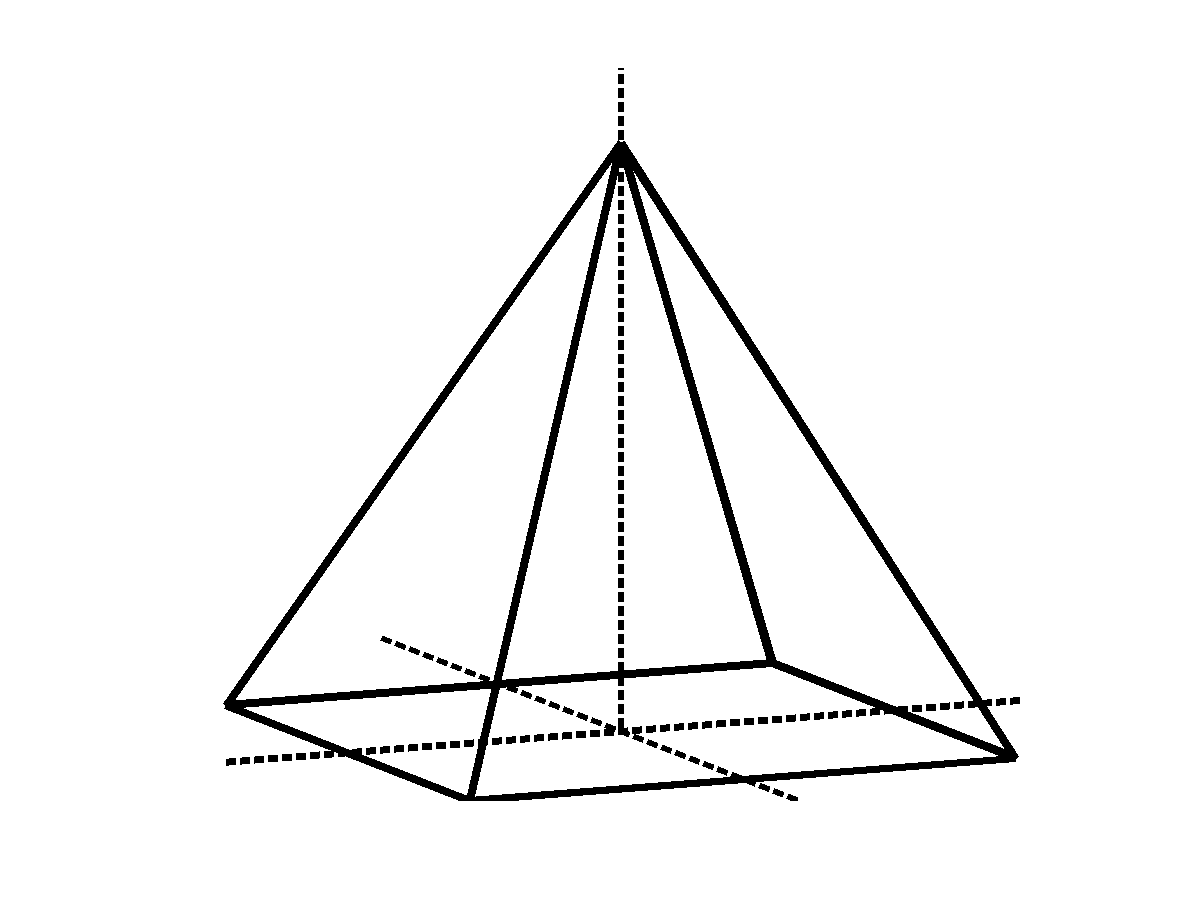}}
\subfigure{\includegraphics[width=.475\textwidth]{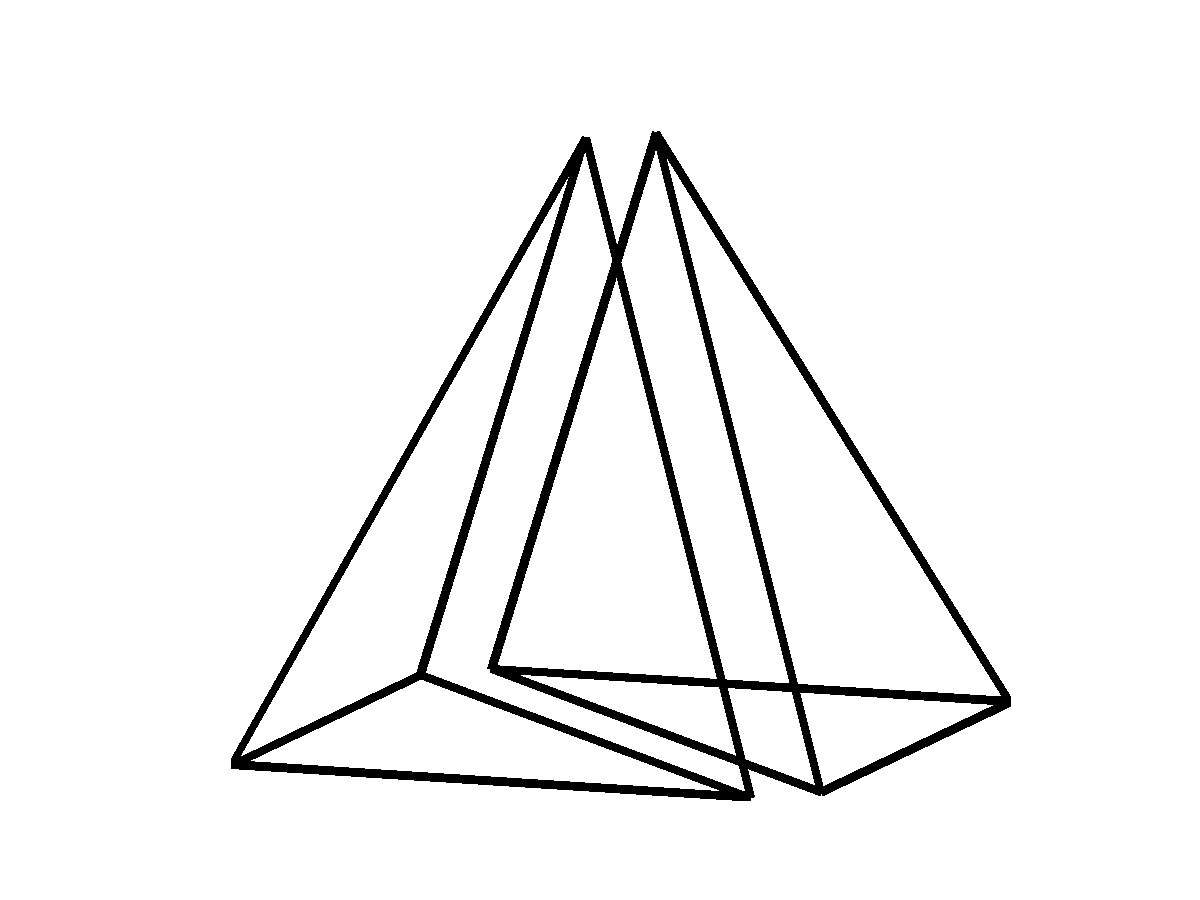}}
\caption{Reference pyramid on $[-1,1]^2 \times [0, 1]$ and tetrahedral splitting. }
\label{fig:split}
\end{figure}

We propose a simple approach of splitting the pyramid into two tetrahedra, as in Figure~\ref{fig:split}, on which we construct nodal distributions using the Warp and Blend procedure.  We will refer to this two-tetrahedron representation of as a ``Duplex pyramid'' (referring both to the tetrahedron as a 3D simplex and the definition of ``duplex'' as having two parts) .  Supposing that the nodes along the shared face are only counted once, this produces ${(N+1)(N+2)(2N+3)}/{6}$ nodes, equal to the dimension $N_p$ of the orthogonal basis.  

The idea of the Duplex pyramid is not new; both Wieners \cite{wieners1997conforming} and Bluck and Walker \cite{bluck2008polynomial} approached the construction of basis functions on the pyramid by dividing the pyramid into two tetrahedral and applying a conformity condition on the interface.  The construction of nodes may be approached the same way --- we will construct nodal sets by dividing the pyramid into two tetrahedra but constrain them to have the same locations on the shared face.  However, unlike the splitting techniques of Wieners and Bluth and Walker, the choice of a nodal distribution is independent of the choice of finite element space, depending only on an orthogonal basis defined over the entire pyramid.  To this end, we choose the orthogonal basis of Bergot, Cohen, and Durufle as described in Section~\ref{sec:bases}.  

\subsubsection{The Warp and Blend procedure}
\label{sec:pyramidWB}
\begin{figure}
\centering
\subfigure{\includegraphics[width=.3\textwidth]{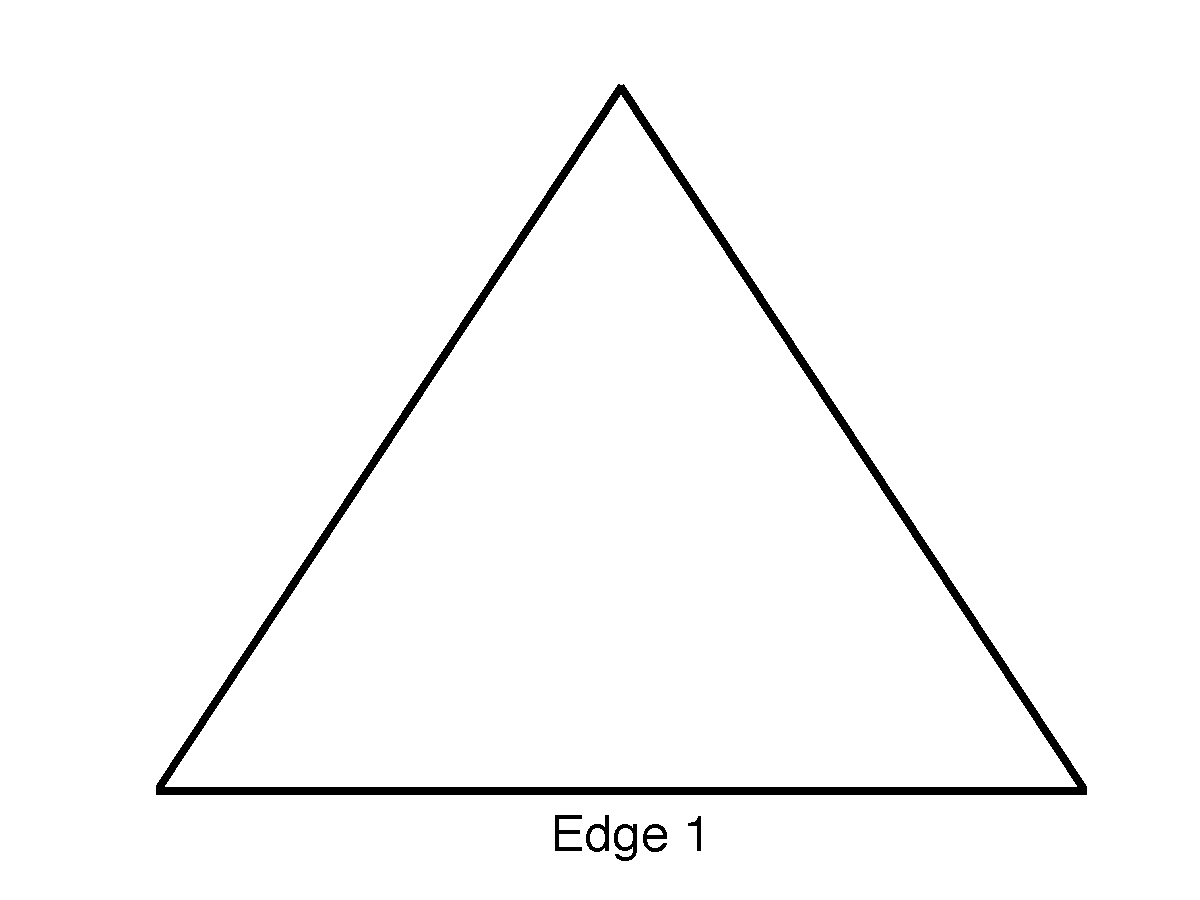}}
\subfigure{\includegraphics[width=.33\textwidth]{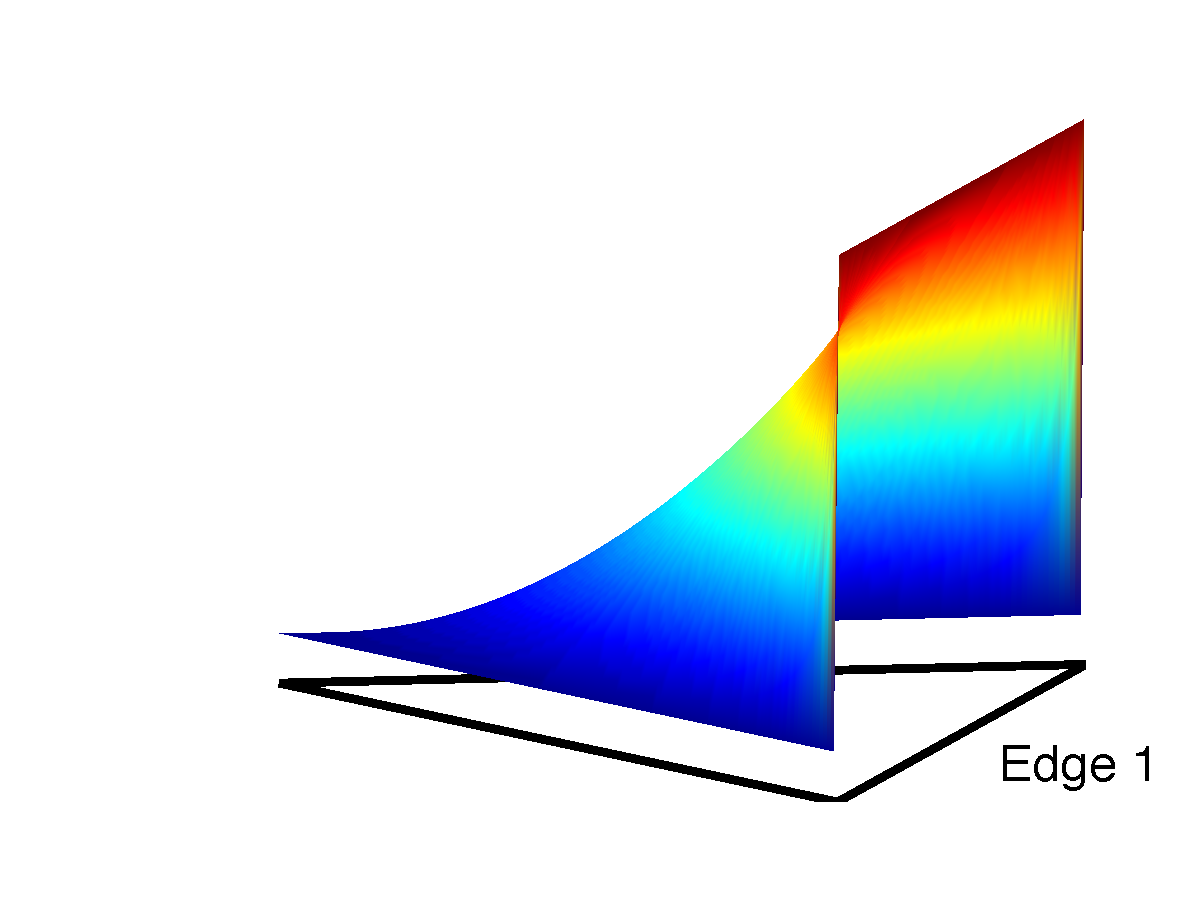}}
\subfigure{\includegraphics[width=.33\textwidth]{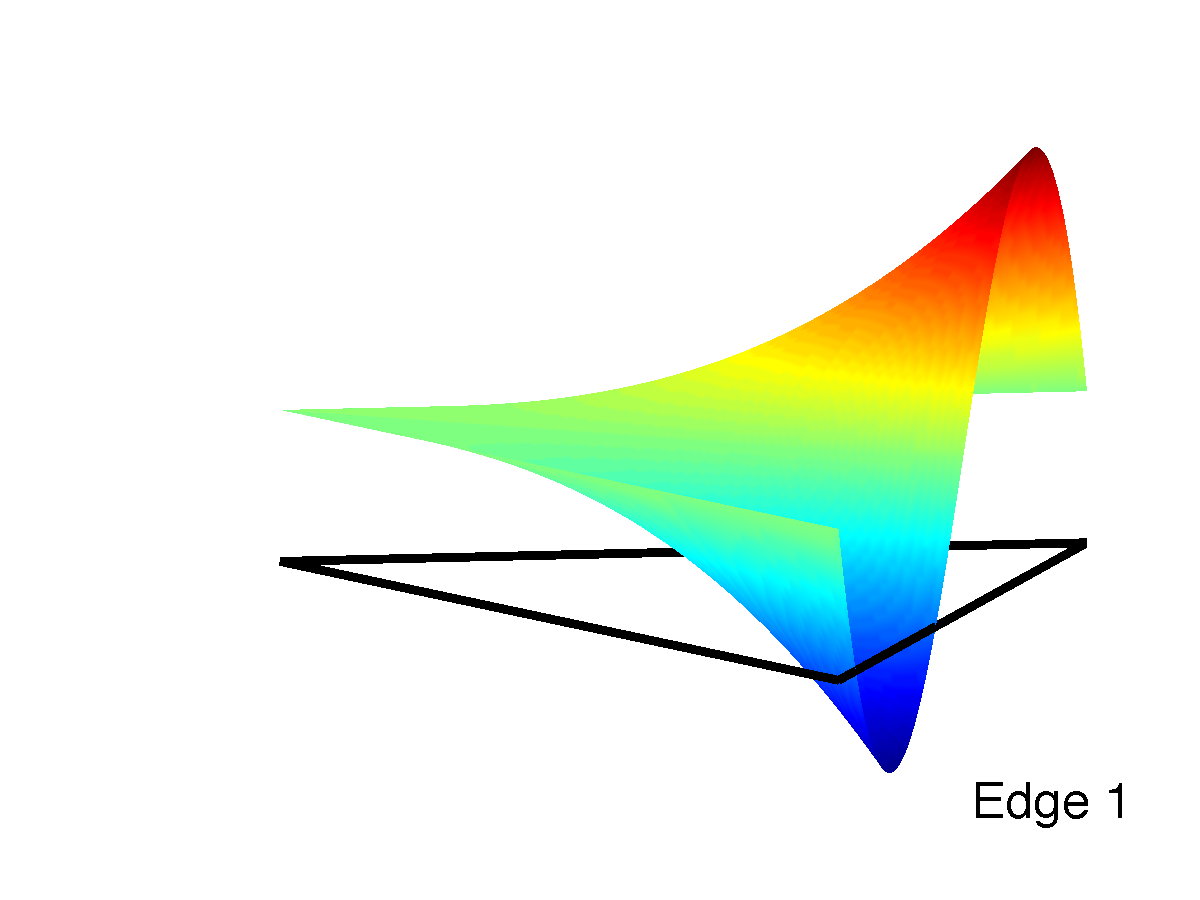}}
\caption{Reference equilateral triangle with $xy$ coordinates (left), blending function, (middle) and combined warp/blend (right) for edge 1.}
\label{fig:tri_blend}
\end{figure}

The Warp and Blend procedure for a triangle is based on the warping function $w_{\rm 1D}(r)$, which maps 1D equidistant points $r_{\rm eq}$ to 1D GLL nodes $r_{\rm GLL}$ via the deformation
\[
r_{\rm GLL} = r_{\rm eq} + w_{\rm 1D}(r_{\rm eq}).
\]
We define $w_{\rm 1D}(r)$ directly as the interpolating polynomial of $r_{\rm GLL} - r_{\rm eq}$.  We may now displace equidistant nodes on a given triangle edge based on the displacement of GLL nodes over an edge.  Consider edge 1; the blend procedure is then to extrapolate the displacement of nodes into the interior of the triangle by defining a blending function of the barycentric coordinates $\lambda_1, \lambda_2, \lambda_3$, with the requirement that the blending function be one on edge 1 and zero on edge 2 and 3.  Figure~\ref{fig:tri_blend} shows the Edge 1 blending function $b_1(\lambda_1,\lambda_2,\lambda_3)$.  

Noting that $1-\lambda_1 -\lambda_2 = \lambda_3$, we may write $b_1$ as
\[
b_1(\lambda_1,\lambda_2,\lambda_3) = \frac{4 \lambda_1 \lambda_2}{(2\lambda_1 + \lambda_3)(2\lambda_2 + \lambda_3)} = \frac{4 \lambda_1 \lambda_2}{1-r_1^2}
\]
where $r_1 = \lambda_2 - \lambda_1 \in [-1,1]$ is the 1D coordinate along edge 1.  

The denominator of $b_1$ is singular at points $r = \pm 1$; however, we note that these correspond to the displacement of vertex nodes by $w_{\rm 1D}(r)$, which we assume is zero for all nodal sets.  This motivates the definition of equivalent nonsingular warping and blending functions $\tilde{b_1}(\lambda_1,\lambda_2,\lambda_3)$ and $\tilde{w}(r)$
\[
\tilde{b}_1 = 4\lambda_1\lambda_2, \qquad \tilde{w}(r) = \begin{cases}
\frac{w_{\rm 1D}(r)}{1-r^2}, & |r| < 1\\
0, & r = \pm 1.
\end{cases}.
\]

The warp is applied along the edge tangent direction, while the blending carries the warp function into the interior of the triangle along the edge normal direction.  For edge 1, the warp thus affects the $x$ position of the triangle nodes, and the blending function carries this $x$ displacement into the interior of the element along the $y$ direction.  This results in the following expression for new nodal positions 
\[
\vectwo{x}{y} = \vectwo{x_{\rm eq}}{y_{\rm eq}} + \vectwo{1}{0}\tilde{w}(\lambda_2-\lambda_1)\tilde{b}_1(\lambda_1,\lambda_2,\lambda_3).
\]
The final step of the Warp and Blend procedure is to parametrize the blending function with a quadratic variation in order to increase the magnitude of the blend in the direction normal to the edge. For edge 1, this may be expressed using the modified blending function $\LRp{1 + (\alpha \lambda_3)^2}\tilde{b}_1$.
As $\alpha$ is increased, the amount which the warping function $\tilde{w}(r)$ is blended towards the opposite vertex increases.  
Taking the same blending function parameter $\alpha$ over all edges allows for a one-parameter family of nodal distributions, which may then be optimized over $\alpha$ to minimize the Lebesgue constant of the resulting nodal distribution.  

For a tetrahedron, since each face is the affine image of an equilateral triangle, we may use the above procedure to define warping and blending functions to displace face nodes.  A face blending function may then be used to define displacement formulas for nodes in the interior of the tetrahedron.  This blend may also be optimized with a some parameter $\beta$, which is arbitrarily taken to be the same as the face blend parameter $\alpha$ to retain the one-parameter nature of the optimization. Since this portion of the Warp and Blend procedure does not change for our extrapolations to the pyramid, we omit the details for brevity and refer the reader to \cite{warburton2006explicit,hesthaven2007nodal}. 

\subsubsection{A Duplex Warp and Blend procedure}
\label{sec:splitTet}
\begin{figure}
\centering
\subfigure{\includegraphics[width=.375\textwidth]{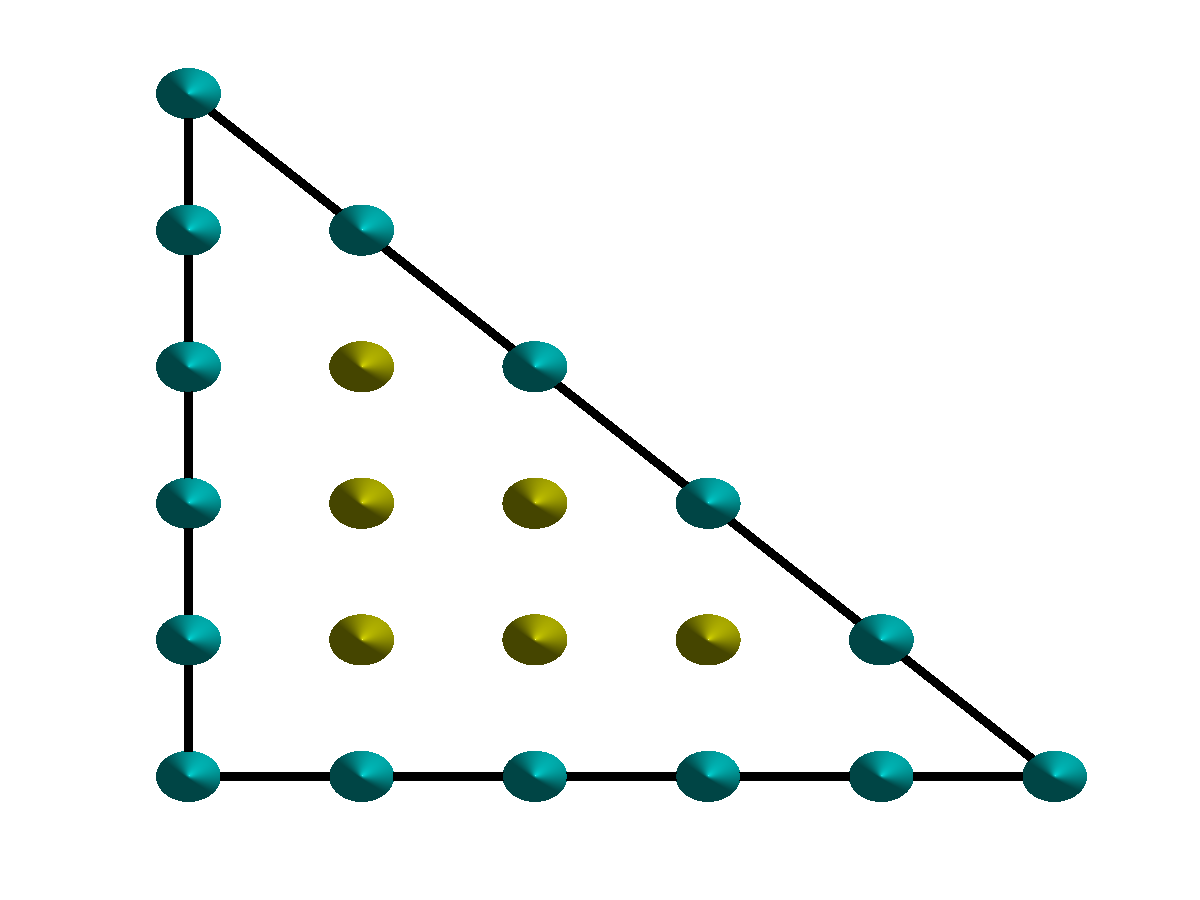}}
\subfigure{\includegraphics[width=.375\textwidth]{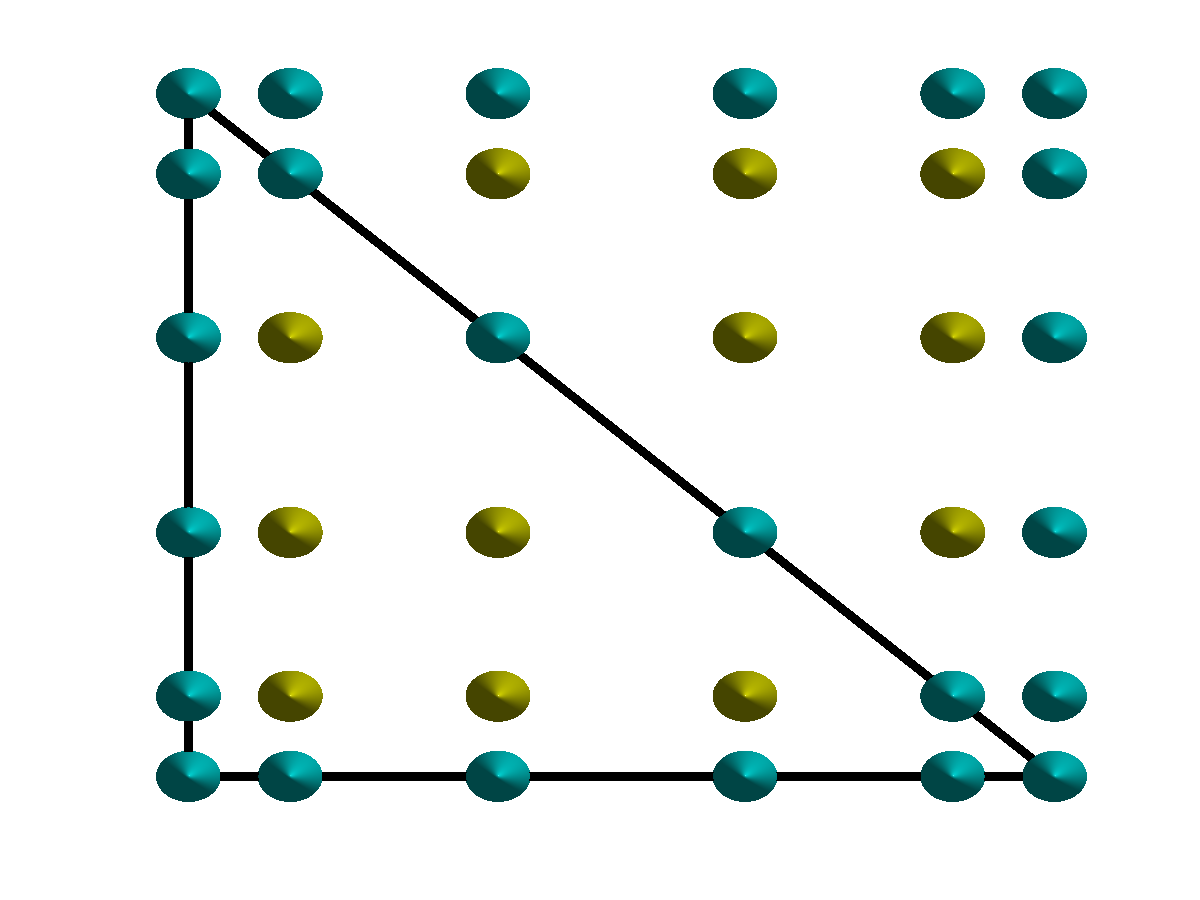}}
\caption{Reference triangle in $rs$ coordinates with equispaced (left) and quadrilateral GLL nodes (right) for $N=5$ overlaid.}
\label{fig:tri_GLL}
\end{figure}

To adapt the Warp and Blend procedure to the Duplex pyramid, we redefine the warping and blending functions for the faces of each tetrahedra corresponding to the square base of the pyramid.  For the $rs$ reference right triangle with coordinates
\[
r,s\in [-1,1] \quad r + s \leq 1
\]
we require the Warp and Blend procedure to map equispaced nodes on a triangle to match half of the GLL nodes on a quadrilateral, as shown in Figure~\ref{fig:tri_GLL}.  We note that, since the GLL nodes are a tensor product in $r$ and $s$ coordinates, we may directly define the face warping function as a tensor product $w_{\rm 1D}(r)w_{\rm 1D}(s)$ as well.  In addition, we assume that the tetrahedral face embedded in the base of the pyramid corresponds to a plane where $t$ is constant, implying that the warp for the base face in the $t$ coordinate is zero.  

Let $x,y,z$ denote coordinates on the equilateral tetrahedron, and let the number of nodes on a tetrahedron be denoted by $N_{p,T}$.  The resulting Warp and Blend procedure for a tetrahedral half of a pyramid is given 
as follows: for faces corresponding to the triangular faces of the pyramid, the Warp and Blend procedure is identical to that of the regular tetrahedron.  For the face corresponding to the square base of the pyramid, we define a warp which maps the nodes to half of a GLL distribution on the right triangle, which is then blended into the interior as usual.  The resulting nodes are then mapped to the two tetrahedral halves of the $[-1,1]^2\times [0, 1]$ reference pyramid using an affine transformation, such that the nodes on the base of each tetrahedra align with tensor product GLL nodes on the square base of the pyramid, and the nodes on the shared face of the two tetrahedra match.  This procedure is given in more detail in Algorithm~\ref{alg:splitWB}.  

\begin{algorithm}
\begin{algorithmic}[1]
\Procedure{Duplex Warp and Blend}{}
\State Initialize $\{x_i,y_i,z_i\}_{i = 1}^{N_{p,T}}$ to equispaced nodes on the tetrahedron.  
\For{faces of the tetrahedron}
\If{face is not the square face}
\State Define face warp $w(x,y,z)$ as for the tetrahedron.
\EndIf
\If{face is the base face}
\State Define $w(x,y,z) = w_{\rm 1D}(r(x,y))w_{\rm 1D}(s(x,y))$.  
\EndIf
\State Blend face warp into interior.
\State Evaluate blended warp, apply shifts to $\{x_i,y_i,z_i\}_{i = 1}^{N_{p,T}}$.  
\EndFor
\State Map $\{x_i,y_i,z_i\}_{i = 1}^{N_{p,T}}$ to $\{r_i,s_i,t_i\}_{i = 1}^{2N_{p,T}}$ on each half of the Duplex pyramid.
\State Remove the $(N+1)(N+2)/2$ redundant nodes on the shared face. 
\State \textbf{return:} points $\{r_i,s_i,t_i\}_{i=1}^{N_p}$.  
\EndProcedure
\end{algorithmic}
\caption{Warp and Blend procedure on one half of the Duplex pyramid.}\label{alg:splitWB}
\end{algorithm}

The Warp and Blend nodes for the triangle and tetrahedron were optimized over a parameter $\alpha$, which controls the quadratic variation of the blending function.  We adopt 
$\alpha = \alpha_{\rm opt}$, the optimized value for the tetrahedron given in \cite{warburton2006explicit, hesthaven2007nodal}.  
\begin{figure}
\centering
\subfigure{\includegraphics[width=.475\textwidth]{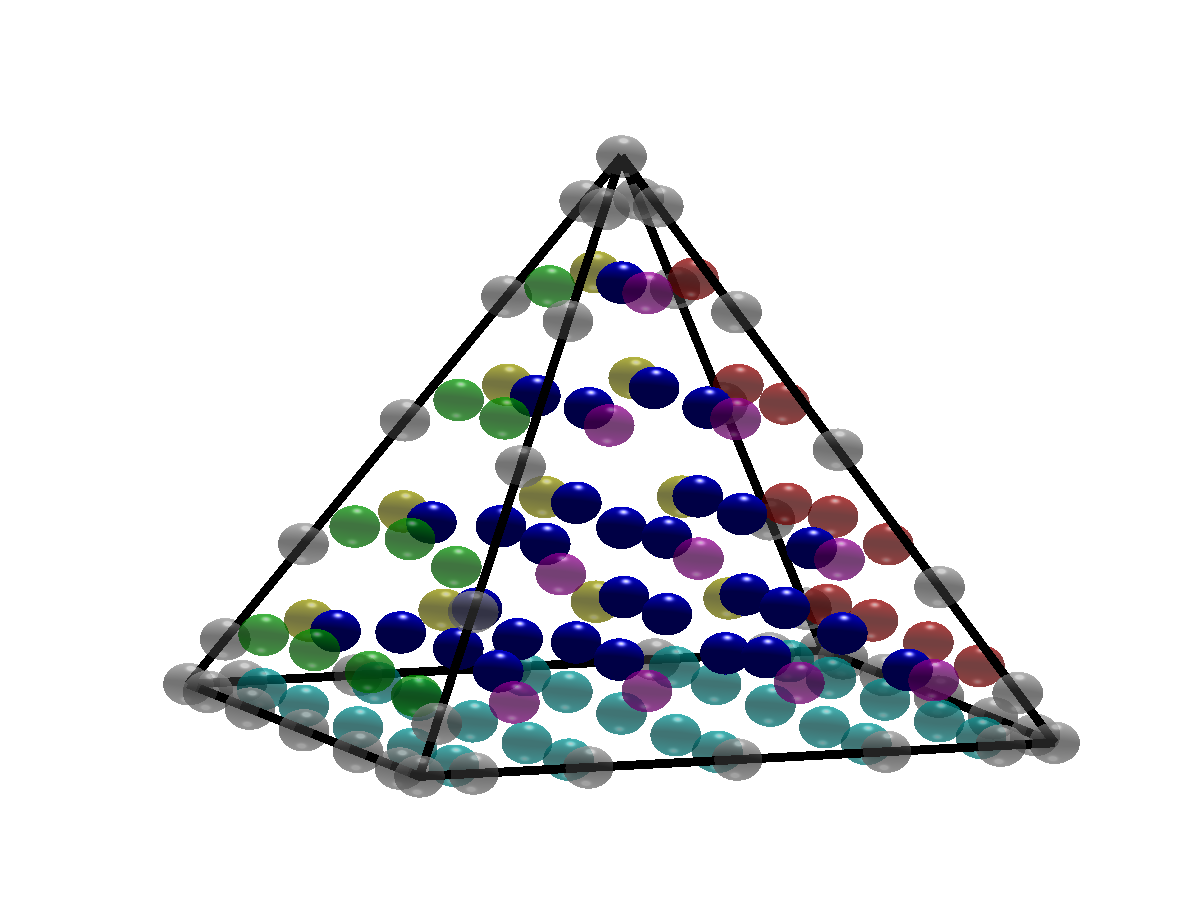}}
\subfigure{\includegraphics[width=.475\textwidth]{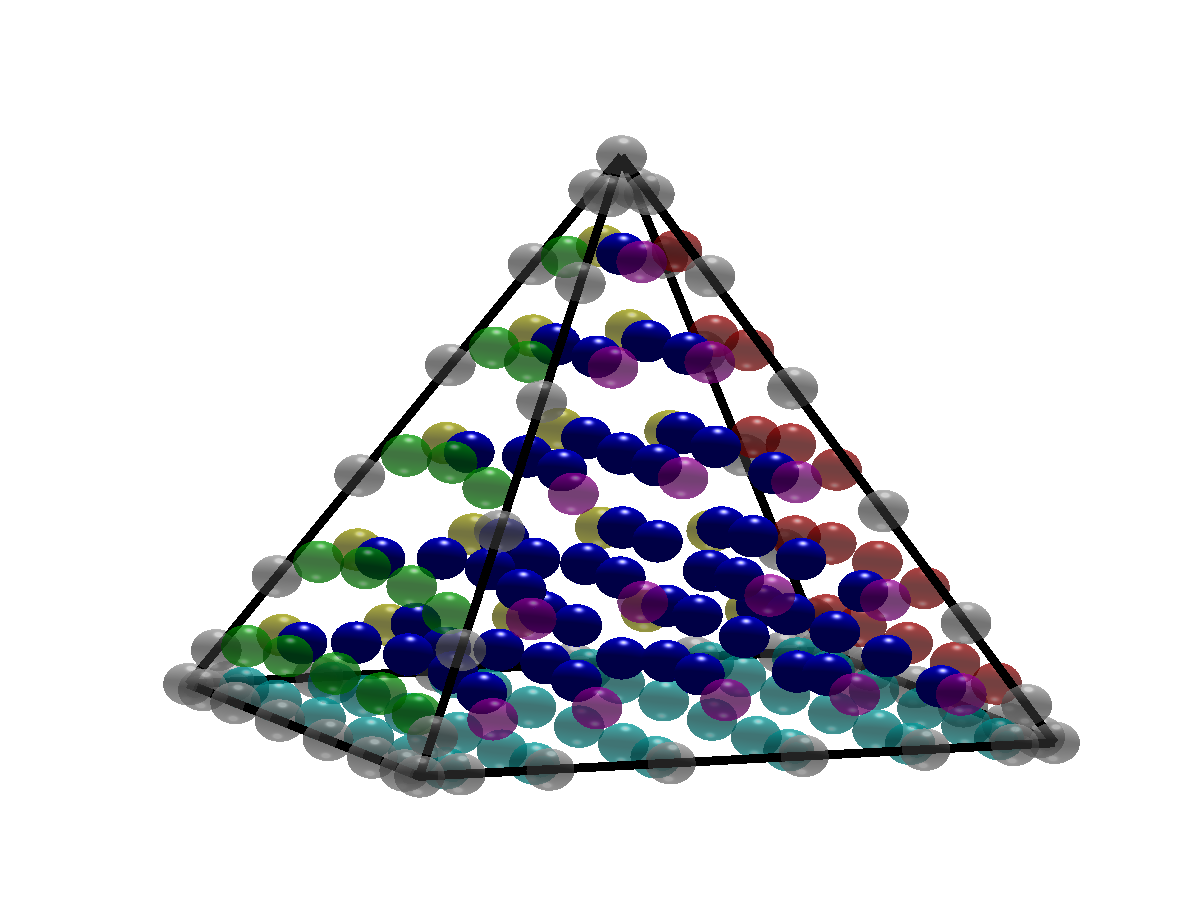}}
\caption{Duplex pyramid Warp and Blend nodes for $N = 6$ (left) and $N=7$ (right) nodes for the pyramid.}
\label{fig:pyrWB}
\end{figure}

The above procedure produces a distribution of nodes for the pyramid which reduces to optimized tetrahedral Warp and Blend nodes on triangular faces and hexahedral GLL nodes on the quadrilateral base, but is not rotationally symmetric in the interior of the pyramid.   We will refer to these nodes as ``Duplex'' in the numerical results.  

\subsection{An Interpolatory Warp and Blend procedure}
\label{sec:InterpWB}
The Duplex pyramid construction of nodes, while a viable procedure for determining a nodal distribution on the pyramid, is less elegant and more complicated than the original Warp and Blend construction on the tetrahedron.  However, it is possible to generalize the Warp and Blend procedure directly to the pyramid in another manner.  We first illustrate the procedure on the triangle.  

In 1D, the warping function $w_{\rm 1D}(r)$ is constructed as the interpolating polynomial of $r_{\rm GLL} - r_{\rm eq}$, and represents a map from equispaced nodes to the difference between the GLL nodes and equispaced nodes.  
The extension to the 2D triangle is discussed in Section~\ref{sec:pyramidWB}.  The warping function for each individual edge is blended into the triangle in order to determine the displacement of nodes in the interior.  The total displacement of the interior nodes may be determined by accumulating the displacements from the warping of each edge.  

In 1D one may also define a direct map from equispaced nodal positions to to GLL nodal positions $r_{\rm GLL} = m_{\rm 1D}(r_{\rm eq})$, where $m_{\rm 1D}(r) = r + w_{\rm 1D}(r)$ is the interpolating polynomial for the positions of the GLL nodes.  We may also define linear vertex shape functions in $r,s$ coordinates 
\[
v_1(r,s) = -\frac{r+s}{2}, \quad, v_2(r,s) = \frac{1+s}{2}, \quad v_3(r,s) = \frac{1+r}{2}.
\]
Since $v_1, v_2, v_3$ are identical to the barycentric coordinates $\lambda_1,\lambda_2, \lambda_3$, we may equivalently use vertex shape functions in lieu of barycentric coordinates in defining a linear blend of $m_{\rm 1D}(r)$ in to the interior of the triangle.  This will be useful in generalizing Warp and Blend to domains without barycentric coordinates.  

Note that for the original Warp and Blend procedure, the blending for each edge is the product of two barycentric variables, which are each linear functions in the local coordinates $r,s$.  As a result, the total displacement of the nodes in each coordinate $r$ and $s$ is the product of linear polynomials and an order $N$ polynomial on the face.  This may be exploited for an Interpolatory Warp and Blend procedure based on the above properties of the map.  

We define the edge basis functions for $j = 0,\ldots, N-2$ using 1D Legendre polynomials $L_j(r)$
\begin{align*}
e_{1,j}(r,s) &= v_1(r,s) v_2(r,s) L_{j}(\xi_1), \quad \xi_1 = v_1(r,s)-v_2(r,s)\\
e_{2,j}(r,s) &= v_2(r,s) v_3(r,s) L_{j}(\xi_2), \quad \xi_2 = v_2(r,s)-v_3(r,s)\\
e_{3,j}(r,s) &= v_3(r,s) v_1(r,s) L_{j}(\xi_3), \quad \xi_3 = v_3(r,s)-v_1(r,s).
\end{align*}
These edge basis functions may be used alongside vertex basis functions and interior (bubble) basis functions that vanish along the boundary to form a hierarchical basis defined over an element, and are commonly used in $hp$-adaptive finite element methods \cite{Demkowicz:2007:CHF:1564840}.  Here, we will discard interior bubble functions, and construct a basis consisting of only vertex and edge shape functions over the surface of the triangle.  

Since the edge basis has cardinality $3N$, equal to the number of nodes on the surface of the triangle, we may use the above basis to interpolate a preset distribution of nodes on the boundary.  The evaluation of these interpolants at equispaced nodes then determines the position of the new Warp and Blend nodes.  

For the triangle, we wish to enforce a GLL nodal distribution over the edge for conformity with quadrilateral elements.  Let $\{\phi_j(r,s)\}_{j = 1}^{3N}$ be the basis consisting of vertex and edge shape functions.  Using the coordinates of the $3N$ GLL points (vertex and edge nodes) over the surface of the triangle, we build a $3N\times 3N$ Vandermonde matrix, which can be used to interpolate, at equispaced points on the triangle surface, the positions of GLL nodes on the edges.  This process explicitly constructs maps from equispaced nodal coordinates to GLL nodal coordinates on each edge, which are blended linearly into the interior.  Interior nodal distributions are then determined as direct evaluations of this map at equispaced coordinates 

\begin{figure}
\centering
\subfigure{\includegraphics[width=.475\textwidth]{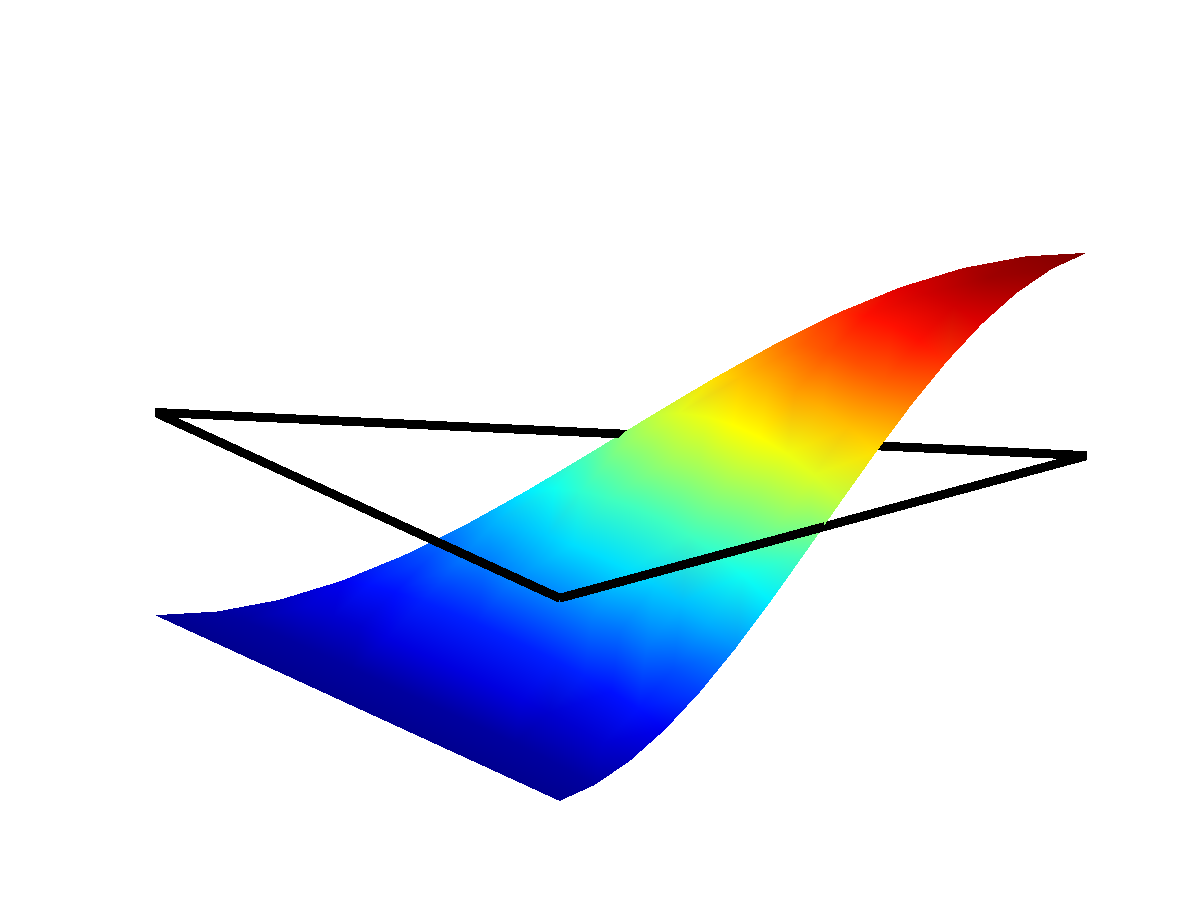}}
\subfigure{\includegraphics[width=.475\textwidth]{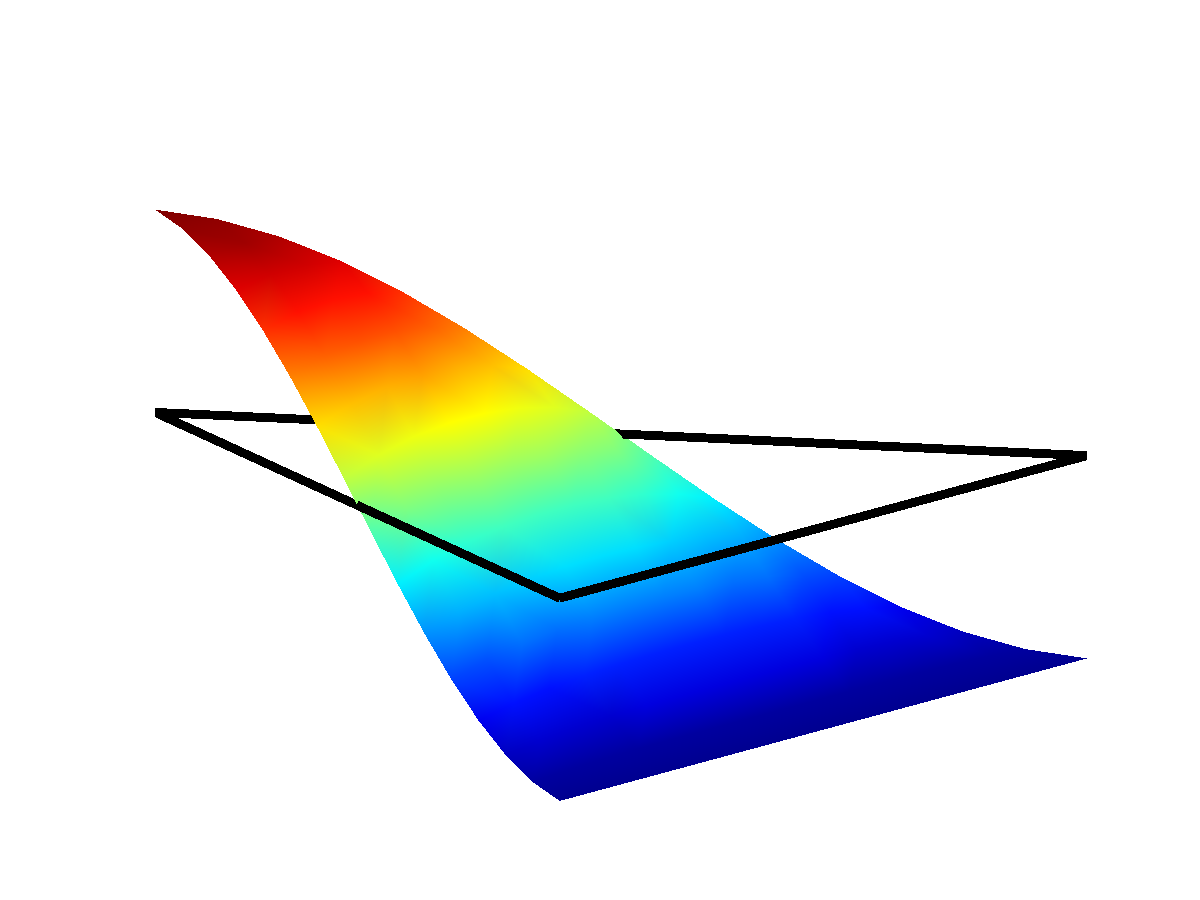}}
\caption{Interpolated maps from equispaced to Warp and Blend nodes for both $r$ (left) and $s$ (right) coordinates on the equilateral triangle.}
\label{fig:wbmaps}
\end{figure}

If we include an additional quadratic blending into the edge basis functions
\begin{align*}
e_{1,j}(r,s) &= \LRp{1 + \LRp{\alpha v_3}^2}v_1(r,s) v_2(r,s) L_{j-1}(\xi_1), \quad \xi_1 = v_1(r,s)-v_2(r,s)
\end{align*}
and similarly for $e_{2,j}, e_{3,j}$, then this process becomes identical to the original Warp and Blend procedure.  In particular, if we set $\alpha = \alpha_{\rm opt}$, the optimized value reported for the triangle in \cite{warburton2006explicit}, we recover exactly (to machine precision) the optimized Warp and Blend nodes on the triangle.


\subsection{An Interpolatory Warp and Blend procedure for the pyramid}

Using the reference pyramid $[-1,1]^2 \times [0,1]$, we may define the vertex shape functions of Bedrosian \cite{bedrosian1992shape}  
\begin{align*}
v_1(r,s,t) &= \frac{1}{4}\left(1-r-s-t+\frac{rs}{1-t}\right), & v_2(r,s,t) &= \frac{1}{4}\left(1+r-s-t-\frac{rs}{1-t}\right)\\
v_3(r,s,t) &= \frac{1}{4}\left(1+r+s-t+\frac{rs}{1-t}\right), &v_4(r,s,t) &= \frac{1}{4}\left(1-r+s-t-\frac{rs}{1-t}\right)\\
v_5(r,s,t) &= t.& 
\end{align*}
Each vertex function vanishes at the other four vertices, and the traces of these shape functions are linear, though the functions themselves are rational.  We may use these vertex functions to generalize the Interpolatory Warp and Blend procedure to the pyramid.  Mimicking the procedure for the triangle, we define a hierarchical basis on both edges and faces of the pyramid.  The edge functions may be defined in a similar manner to the triangle; for an edge between vertices $a$ and $b$, we may define $N-1$ edge functions
\[
e_{ab,j}(r,s,t) = v_a(r,s,t) v_b(r,s,t) L_{j}(\xi_{ab}), \quad j = 0,\ldots, N-2,
\]
where $L_j(\xi)$ is again the $j$th order 1D Legendre basis function and $\xi_{ab} = v_a(r,s,t) - v_b(r,s,t)$ is the local coordinate along the edge.  

We may similarly define triangular face functions in terms of vertex shape functions due to the linearity of their traces.  For a triangular face defined by vertices $a$, $b$, and $c$, we may define barycentric coordinates over the face in terms of vertex shape functions
\[
\lambda_1(r,s,t) = v_a(r,s,t), \quad \lambda_2(r,s,t) = v_b(r,s,t), \quad \lambda_3(r,s,t) = v_c(r,s,t).
\]
These may then be used to evaluate the orthogonal Dubiner basis on the triangle.  Let $D_{j}(\lambda_1,\lambda_2,\lambda_3)$ denote the $j$th Dubiner polynomial  as a function of the barycentric coordinates; for $j = 0,\ldots, {(N - 1)(N - 2)}/{2}$, we may then define the triangular face functions through 
\[
f_{{\rm tri},j}(r,s,t) = v_a(r,s,t) v_b(r,s,t) v_c(r,s,t) D_{j}(v_1,v_2,v_3).
\]

For $0 \leq j, k \leq N-1$, the pyramid base face functions may also be defined using the vertex functions at the base and tensor products of 1D Legendre polynomials
\[
f_{{\rm quad},jk}(r,s,t) = v_a(r,s,t) v_b(r,s,t) v_c(r,s,t) v_d(r,s,t) L_j(r) L_k(s).
\]
Since both the number of surface nodes and total number of vertex, edge, and face basis functions are $N_p = 3N^2 + 2$, we may define the basis $\{\phi_j(r,s,t)\}_{j=1}^{N_p}$ as the collection of vertex, edge, and face functions and construct a square Vandermonde matrix $V$ over the surface nodes.\footnote{We have constructed the above basis to closely mimic the blending functions used in the original Warp and Blend procedure.  However, we note that any choice of hierarchical basis defined over vertices, edges, and faces --- for example, the $H^1$-conforming basis described in Nigam and Phillips \cite{nigam2012high}, Bergot, Cohen, and Durufle \cite{bergot2010higher}, or \cite{karniadakis2013spectral, Demkowicz:2007:CHF:1564840} --- could also be used to construct maps from equispaced to Warp and Blend nodes.}

The construction of the map from equispaced nodes to Warp and Blend nodes on the face may then be expressed in coefficients of the basis $\{\phi_j(r,s,t)\}$.  
Assuming that the pyramid surface nodes are a combination of Warp and Blend nodes on triangular faces and tensor product GLL nodes on the quadrilateral faces, we may solve for the interpolant of these surface nodal values at equispaced nodes on the surface of the pyramid.  These interpolants are regarded as maps whose evaluation at equispaced points on the surface and interior of the pyramid determine the position of the Interpolatory Warp and Blend nodes.  The procedure for computing these nodes is given in Algorithm~\ref{alg:IWB}, and we will refer to these nodes as ``IWB'' in the numerical results.  
\begin{algorithm}
\begin{algorithmic}[1]
\Procedure{Interpolatory Warp and Blend}{}
\State Initialize positions of  equispaced nodes on the pyramid $r_i^{\rm eq},s_i^{\rm eq},t_i^{\rm eq}$ .
\State Select $3N^2+2$ target node positions on the faces of the tetrahedron.  
\State Define the vertex, edge, and face basis functions $\{\phi_j(r,s,t)\}_{j=1}^{3N^2 + 2}$.  
\State Using $\phi_j$, compute interpolants $m_r(r,s,t), m_s(r,s,t), m_t(r,s,t)$ of the target nodal positions.  
\State Evaluate $r_i = m_r(r_i^{\rm eq}, s_i^{\rm eq}, t_i^{\rm eq}), \quad s_i = m_s(r_i^{\rm eq}, s_i^{\rm eq}, t_i^{\rm eq}), \quad t_i = m_t(r_i^{\rm eq}, s_i^{\rm eq}, t_i^{\rm eq}).$
\State \textbf{return:} points $\{r_i,s_i,t_i\}_{i=1}^{N_p}$.  
\EndProcedure
\end{algorithmic}
\caption{Interpolatory Warp and Blend procedure for the pyramid.}
\label{alg:IWB}
\end{algorithm}

\begin{figure}
\centering
\subfigure{\includegraphics[width=.475\textwidth]{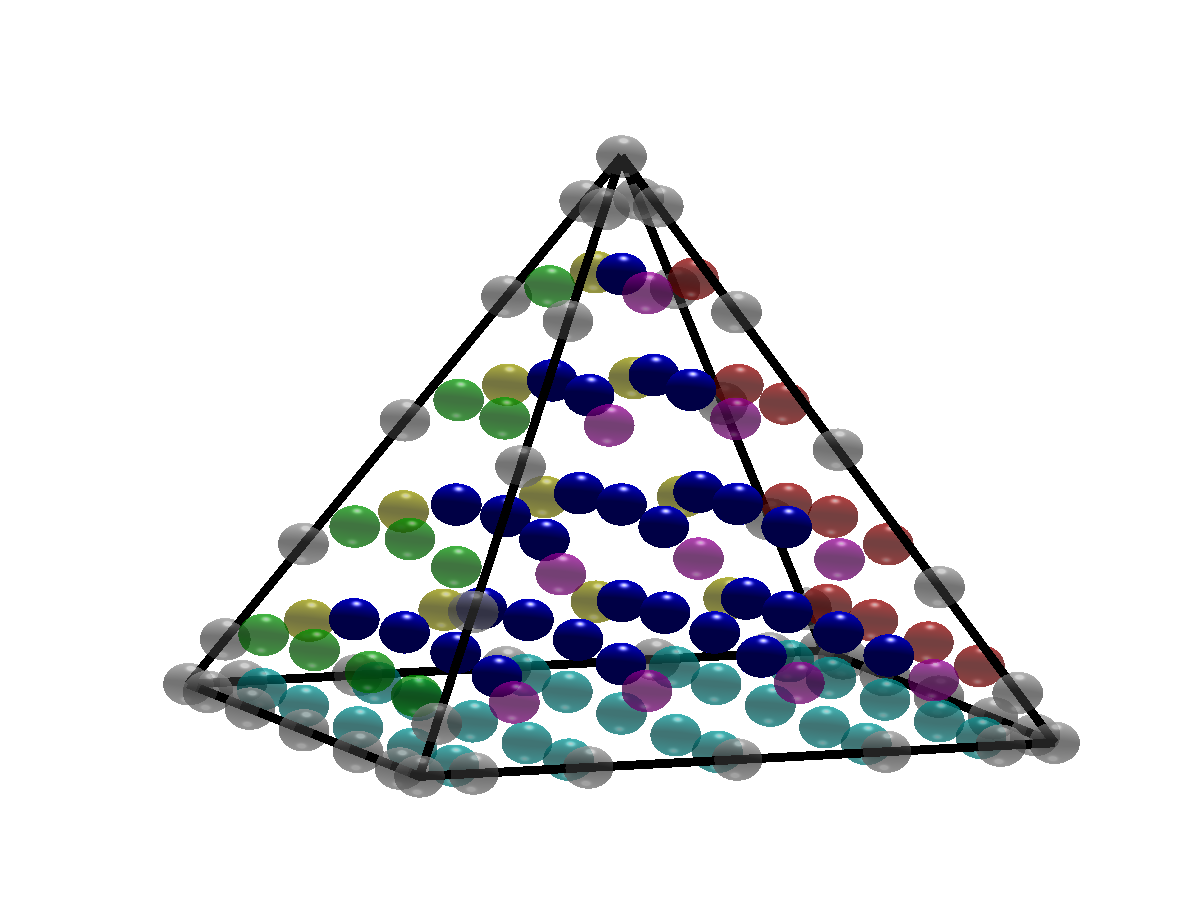}}
\subfigure{\includegraphics[width=.475\textwidth]{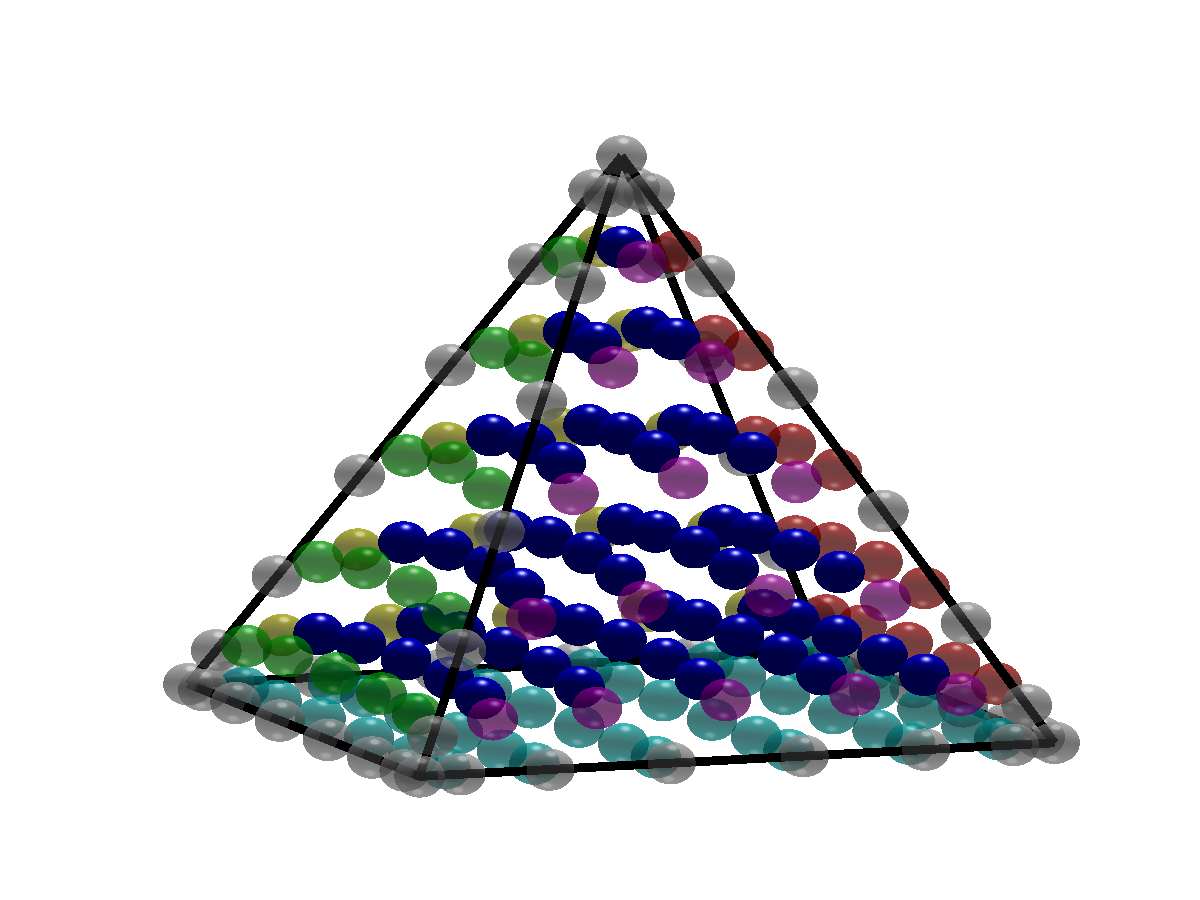}}
\caption{Interpolatory Warp and Blend nodes for $N = 6$ (left) and $N=7$ (right) nodes for the pyramid.}
\label{fig:pyrIWB}
\end{figure}

While this interpolatory procedure is identical to the original Warp and Blend procedure in 2D, they contain differences in 3D.  The original Warp and Blend procedure modified node positions edge by edge: for a given triangular face, the warp from each edge of the triangle is blended into the face.  Since this warp represents the displacement needed to move an equispaced node to a Warp and Blend node, these warps may be applied edge-by-edge, updating the positions of nodes on the face one edge at a time.  For a given edge on a given face, this also defines a face warp, which is then blended into the interior via a face blending function which vanishes on all other faces.  Conceptually, the warping of interior nodes is related to the warping of edge nodes only indirectly (through the face warp) --- edge and interior nodes are decoupled from each other.  

In contrast, in the Interpolatory Warp and Blend procedure, the warping of edge nodes directly affects the interior node distribution, due to the fact that each edge function is blended into the using one or two vertex shape functions, which are also nonzero in the interior of the pyramid.  Thus, both face and edge nodes are coupled together in determining the position of interior nodes of the pyramid.

While the Fekete nodes tend to perform deliver the lowest Lebesgue constant for large $N$, they become more computationally challenging to determine as $N$ increases.  In contrast, apart from an optional 1D optimization of the nodal distribution, both the Duplex and Interpolatory Warp and Blend procedures give explicit (non-iterative) constructions of pyramid nodes for any degree $N$.  

\section{Numerical experiments and comparisons with existing nodal sets}

In this section, we compare Fekete, ``approximate Fekete'', and Warp and Blend nodes with other unoptimized nodal sets for the pyramid.  We examine four metrics: the Lebesgue constant, the determinant of the Vandermonde matrix, the condition number of the Vandermonde matrix, and interpolation error for two specific functions.  In all cases, the Vandermonde matrix is normalized by the $L^2$ norm of the corresponding basis function, computed using quadrature \cite{stroudCode}.  

Since there is no closed form expression through which to explicitly compute the Lebesgue constant for a nodal set, we adaptively sample the Lebesgue function using a random search \cite{warburton2006explicit}
\[
L(\vect{x}) = \sum_{i = 1}^{N_p} \left|\ell_i(\vect{x})\right|
\]
and seek the Lebesgue constant as the maximum value $\Lambda = \max_{\vect{x}\in K} L(\vect{x})$.  

\begin{figure}
\centering
\subfigure{\includegraphics[width=.32\textwidth]{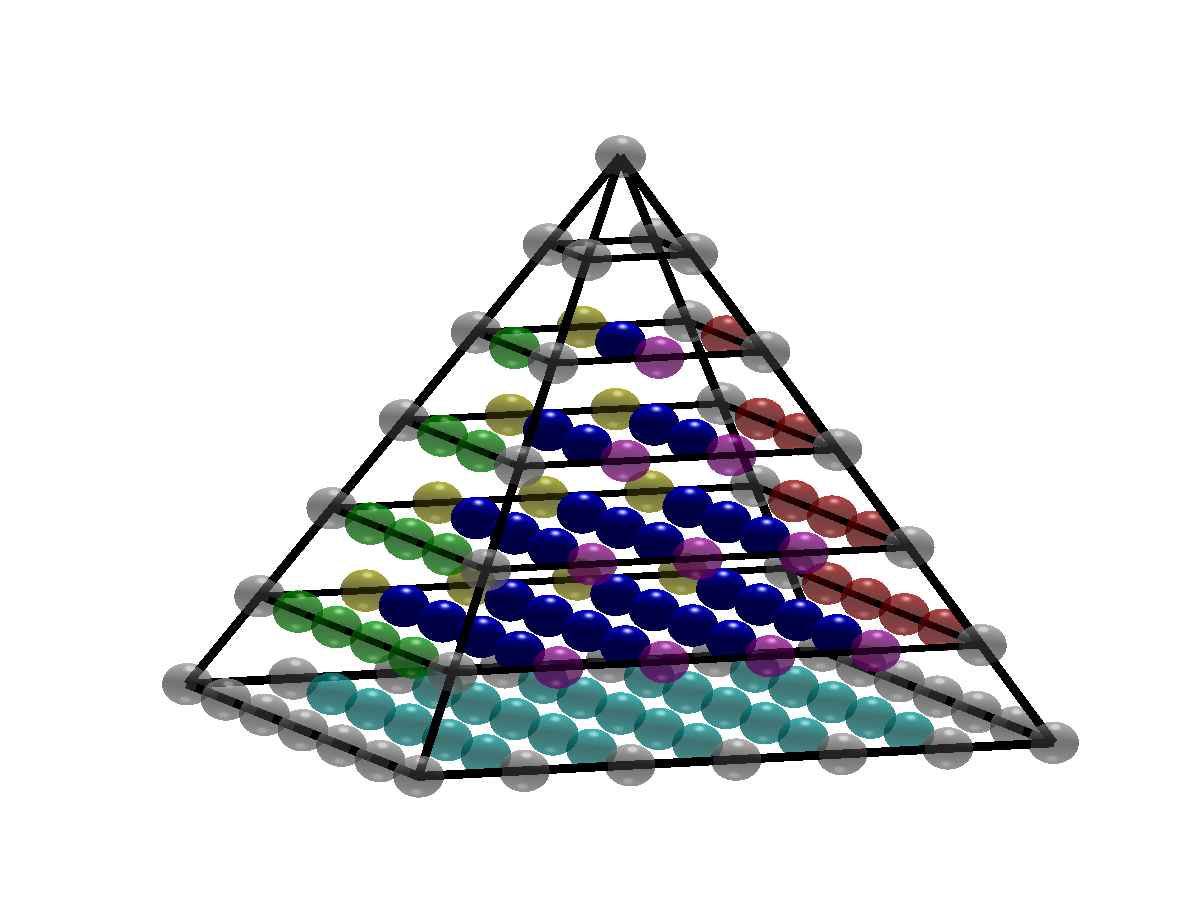}}
\subfigure{\includegraphics[width=.32\textwidth]{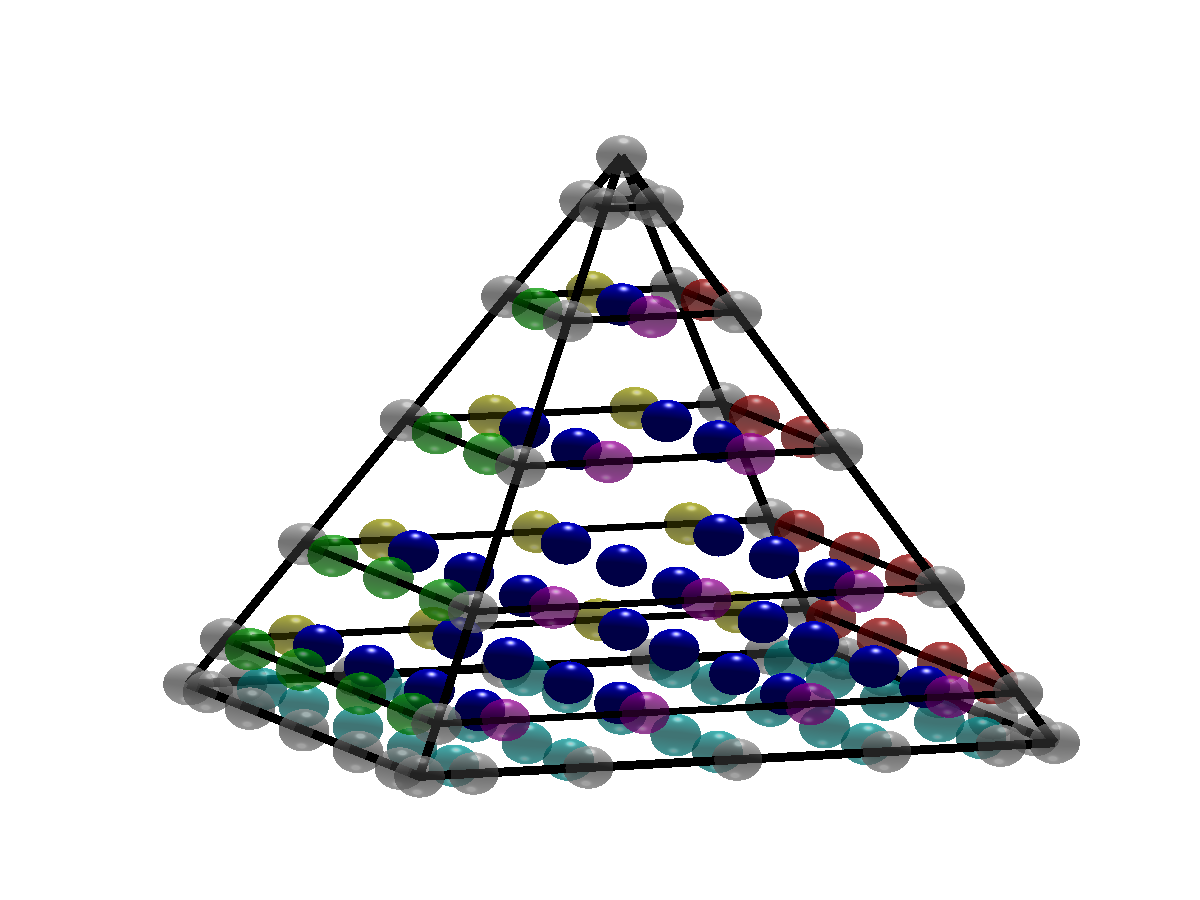}}
\subfigure{\includegraphics[width=.32\textwidth]{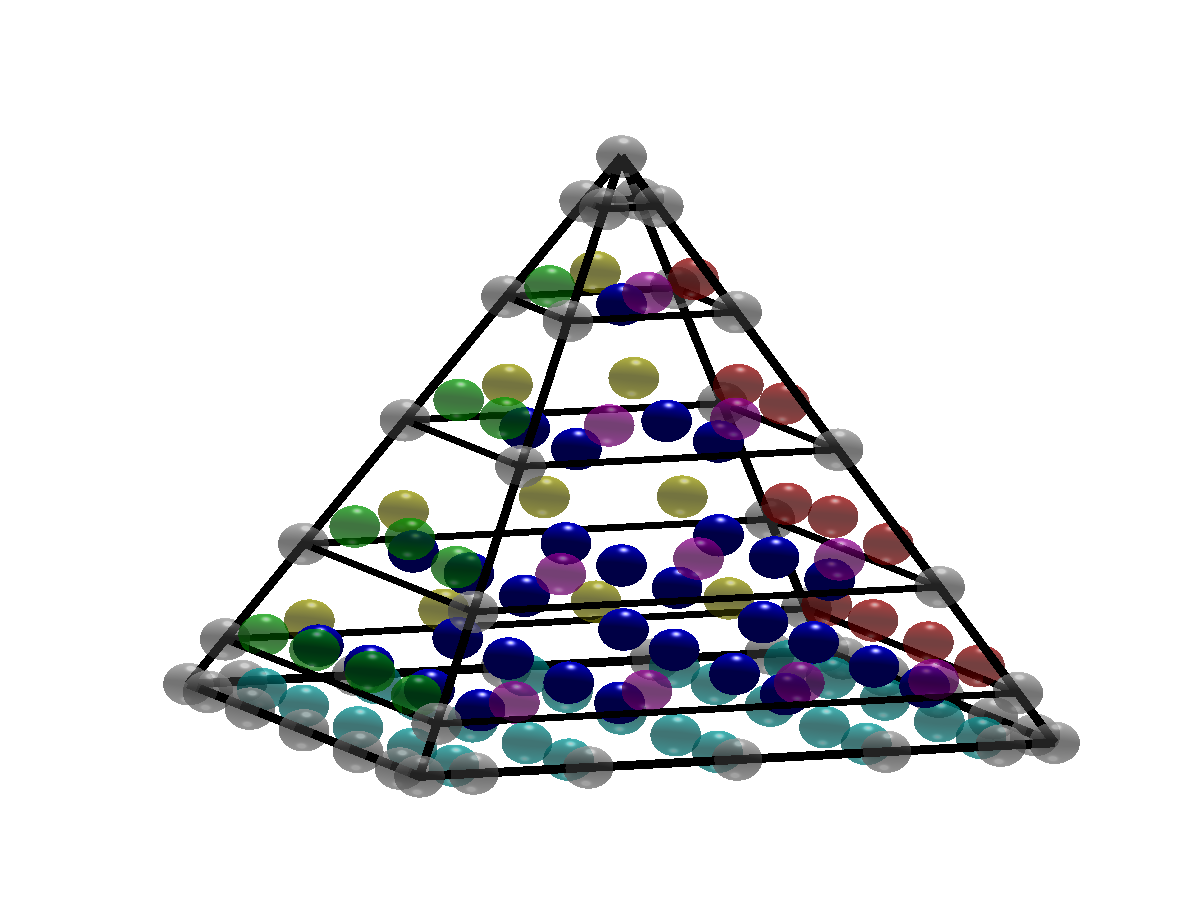}}
\caption{Equispaced (left), Stroud-type GLL (middle), and mixed Warp and Blend/GLL (right) pyramid nodes for $N=6$.}
\label{fig:stroud}
\end{figure}

The baseline comparison is with equispaced nodes on the pyramid; these can be defined level-by-level, similar to the manner in which Stroud conical quadrature rules are constructed \cite{hammer1956numerical}.  Similarly, we may construct a Stroud-type conical GLL node set by levels as well --- the levels are placed according to a GLL distribution, and on each level, nodes are arranged as a tensor product of GLL nodes.  Both are shown in Figure~\ref{fig:stroud} for $N=6$, and are referred to as ``Equi'' and ``Conical'' in the numerical results.    

Both Bergot et al \cite{bergot2010higher} and Gassner et al \cite{gassner2009polymorphic} used electrostatic nodes on the faces and a Stroud-type GLL nodal distribution in the interior (shown in Figure~\ref{fig:stroud}).  For consistency in comparison with our optimized nodal sets, we will mimic their choice of GLL interior nodes but substitute Warp and Blend nodes from the tetrahedron for the electrostatic nodes on the faces.  We refer to this distribution as ``Face'' in the numerical results.  

\subsection{Lebesgue constants}

We include Table~\ref{table:legend}, which describes references and summaries of the various surveyed nodal sets.  
\begin{table}
\centering                                                                                   
\begin{tabular}{| l | l |}
\hline
Equi & Equispaced nodes (see Fig~\ref{fig:stroud})\\
\hline
Conical & Stroud-type (see Fig~\ref{fig:stroud})\\
\hline
Face & Tet face nodes, GLL interior nodes (see Fig~\ref{fig:stroud}, \cite{bergot2010higher,gassner2009polymorphic})\\ 
\hline
Fekete & Steepest ascent ODE-based (see Fig~\ref{fig:pyrFek}, \cite{taylor2000algorithm})\\ 
\hline
Greedy & Choosing from sampled points  \cite{Sommariva20091324, Briani20122477}\\ 
\hline
QR & Iterative refinement using QR  \cite{Sommariva20091324, Briani20122477}\\ 
\hline
Duplex & Two-tet Warp and Blend(see Section~\ref{sec:splitTet})\\ 
\hline
IWB & Interpolatory Warp and Blend (see Section~\ref{sec:InterpWB})\\
\hline
\end{tabular}
\caption{Legend of abbreviations and summaries for different nodal sets.}
\label{table:legend}
\end{table}
Table~\ref{table:leb} gives Lebesgue constants for  different nodal distributions with $N= 3,\ldots,10$.  Since the nodal distribution for $N = 1,2$ is the same for all distributions, we focus on the lowest Lebesgue constants for $N>2$, which are bolded for reference.  

For $N=3$, the Face nodes return the lowest Lebesgue constant by $.01$.  The Duplex pyramid nodes contain the lowest Lebesgue constants for $3 < N \leq 6$, while similarly to the triangle, Fekete nodes outperform other nodal sets for high $N$.  While this happens for $N>10$ on the triangle, it occurs earlier at $N>7$ on the pyramid.  

An interesting observation is that the Face distribution, which chooses face and interior nodal distributions independently of each other, results in a Lebesgue constant which grows faster for large $N$ than the Lebesgue constant for the Conical GLL distribution, implying that the interior distribution of nodes plays a significant role in minimizing the Lebesgue constant.  

\begin{table}[h!]                                                     
\centering                                                            
\begin{tabular}{|c|c|c|c|c|c|c|c|c|}                                  
\hline               
N & Equi & Conical & Face & Fekete & Greedy & QR & Duplex & IWB  \\ 
\hhline{|=|=|=|=|=|=|=|=|=|}                                               
3 & 3.15 & 2.83 & $\bf 2.72$ & 2.73 & 2.80 & 2.80 & 2.73 & 2.75 \\          
\hline                                                                
4 & 5.94 & 4.29 & 4.22 & 4.13 & 4.19 & 4.19 & $\bf 3.80$ & 3.90 \\          
\hline                                                                
5 & 11.87 & 6.84 & 6.93 & 5.53 & 6.33 & 6.03 & $\bf 5.06$ & 5.11 \\         
\hline                                                                
6 & 25.13 & 10.10 & 10.67 & 7.35 & 8.51 & 8.29 & $\bf 6.66$ & 7.23 \\       
\hline                                                                
7 & 56.66 & 14.20 & 15.42 & 9.71 & 12.82 & 13.63 & $\bf 9.65$ & 9.75 \\     
\hline                                                                
8 & 136.40 & 20.43 & 22.16 & $\bf 12.79$ & 18.85 & 21.43 & 14.65 & 14.22 \\ 
\hline                                                                
9 & 350.23 & 31.14 & 34.28 & $\bf 17.16$ & 22.84 & 33.31 & 23.39 & 20.82 \\ 
\hline                                                                
10 & 954.08 & 48.38 & 54.27 & $\bf 25.50$ & 42.85 & 37.23 & 40.12 & 32.16 \\
\hline                                                  
\end{tabular}                                                         
\caption{Values of the Lebesgue constant $\Lambda$ for various nodal sets and polynomial orders $N$.  The smallest Lebesgue constants for $N\geq 3$ are bolded. }
\label{table:leb}
\end{table}   

\subsection{Determinant of the Vandermonde matrix}

In this section, we discuss the magnitude of determinants of the Vandermonde matrix for various nodal sets.  Since the determinants of the normalized Vandermonde matrix are too large to represent numerically, we arbitrarily scaled the Vandermonde matrix 
to prevent numerical overflow.  
As expected, the largest magnitude determinants of the Vandermonde matrix are produced by the Fekete points.  The ``approximate Fekete'' points do remarkably well, producing determinants that are within a factor of the magnitude of the Fekete determinant.  The determinants of all other sets behave roughly the same, decreasing at a steady rate as $N$ increases.

\subsection{Conditioning of the Vandermonde matrix}

In this section, we compare the condition numbers of the normalized Vandermonde matrix for different nodal sets at various $N$.  Since nodal basis functions are typically constructed through the inversion of a Vandermonde matrix, poor conditioning can result in the loss of accuracy when constructing interpolants.  

Overall, the Duplex nodes produce the most well-conditioned matrices, followed closely by the two ``approximate Fekete'' node sets.  The Fekete and Interpolatory Warp and Blend nodes surprisingly do the most poorly apart from Equispaced nodes; however, even in these cases, the condition number is relatively small and should not introduce numerical issues.  

\begin{table}
\centering                                                                                    
\begin{tabular}{|c||c|c|c|c|c|c|c|c|}                                                          
\hline
N & Equi & Conical & Face & Fekete & Greedy & QR & Duplex & IWB  \\ 
\hhline{|=|=|=|=|=|=|=|=|=|}
\hline                                                                        
3 & 15.84 & 16.43 & 16.34 & 16.15 & 16.96 & 16.96 & 16.15 & 16.37 \\          
\hline                                                                        
4 & 22.15 & 20.57 & 20.57 & 20.82 & 21.91 & 21.91 & 20.10 & 20.63 \\          
\hline                                                                        
5 & 34.79 & 29.69 & 29.64 & 27.85 & 34.41 & 31.82 & 26.53 & 28.74 \\          
\hline                                                                        
6 & 60.84 & 38.26 & 38.48 & 40.61 & 48.40 & 43.17 & 36.69 & 43.53 \\          
\hline                                                                        
7 & 123.46 & 53.02 & 53.77 & 63.91 & 51.75 & 76.87 & 52.20 & 67.54 \\         
\hline                                                                        
8 & 301.65 & 80.85 & 83.67 & 107.05 & 82.14 & 114.65 & 78.92 & 110.58 \\      
\hline                                                                        
9 & 810.06 & 131.13 & 135.18 & 188.56 & 105.39 & 191.26 & 124.46 & 187.81 \\  
\hline                                                                        
10 & 2346.19 & 222.97 & 234.38 & 345.23 & 226.88 & 205.85 & 202.36 & 330.31\\
\hline     
\end{tabular}                                                                                 
\caption{Condition numbers of the normalized Vandermonde matrix for various nodal sets and polynomial orders $N$.}
\label{table:VDMcond}
\end{table}                 
         
\subsection{Interpolation errors}

In this section, we compute interpolation errors in the max norm for two functions using our construced nodal sets.  We iterate towards the max norm error $\|f-f_N\|$ using an adaptive sampling, similar to the manner in which the Lebesgue constant is computed.  Tables~\ref{table:interp1} and \ref{table:interp2} shows interpolation errors for two functions: a smooth analytic function $f_1$ and a Runge-type function $f_2$
\begin{align}
f_1(r,s,t) &= (r+1)(s+1)(t+1)\cosh(r+s+t - 1) \label{f1.eqn}\\
f_2(r,s,t) &=\frac{1}{1 + (r^2 + s^2 + t^2)/2}.\label{f2.eqn}
\end{align} 
    
\begin{table}[h!]                                                                                    
\centering                                                                                           
\begin{tabular}{|c||c|c|c|c|c|c|c|c|}                                                                 
\hline                                                                                               
N & Equi & Conical & Face & Fekete & Greedy & QR & Duplex & IWB   \\ 
\hhline{|=|=|=|=|=|=|=|=|=|}                                               
4 & 3.2e-2 & 2.17e-2 & 2.13e-2 & 2.11e-2 & 2.11e-2 & 2.17e-2 & 2.14e-2 & 2.20-2 \\ 
\hline                                                                                               
5 & 8.6e-3 & 4.95e-3 & 4.93e-3 & 5.03e-3 & 4.86e-3 & 4.98e-3 & 4.64e-3 & 4.70e-3 \\ 
\hline                                                                                               
6 & 1.5e-3 & 6.55e-4 & 6.46e-4 & 6.53e-4 & 6.25e-4 & 6.05e-4 & 6.12e-4 & 6.51e-4 \\ 
\hline                                                                                               
7 & 2.2e-4 & 8.27e-5 & 6.99e-5 & 7.64e-5 & 7.16e-5 & 7.12e-5 & 7.56e-5 & 6.93e-5 \\ 
\hline                                                                                               
8 & 3.5e-5 & 1.27e-5 & 1.04e-5 & 1.05e-5 & 1.06e-5 & 1.02e-5 & 9.61e-6 & 1.05e-5 \\ 
\hline                                                                                               
9 & 3.1e-6 & 8.32e-7 & 6.82e-7 & 4.55e-7 & 6.85e-7 & 6.65e-7 & 6.45e-7 & 5.99e-7 \\ 
\hline                                                                                               
10 & 4.8e-7 & 1.69e-7 & 1.18e-7 & 9.23e-8 & 1.31e-7 & 1.31e-7 & 1.27e-7 & 1.26e-7 \\
\hline                                                                                               
\end{tabular}                                                                                        
\caption{Interpolation errors for $f_1$ (see Equation \ref{f1.eqn}).}
\label{table:interp1}                                                             
\end{table}               

We can observe from the results in Table~\ref{table:interp1} that for $N=10$, the Fekete nodes give back the lowest interpolation error for the smooth analytic function, which is consistent with the Fekete nodes having the lowest Lebesgue constant at high $N$.  Equispaced and Conical/GLL nodes deliver larger interpolation errors, especially as $N$ increases.  However, the Face node set (Warp and Blend faces with GLL interior nodes) do quite well for even $N$, which may be specific to the specific function $f_1$.  
                                          
\begin{table}[h!]                                                                                    
\centering                                                                                           
\begin{tabular}{|c||c|c|c|c|c|c|c|c|}                                                                 
\hline                                                                                               
N & Equi & Conical & Face & Fekete & Greedy & QR & Duplex & IWB  \\ 
\hhline{|=|=|=|=|=|=|=|=|=|}                                               
4 & 5.2e-3 & 5.95e-3 & 6.04e-3 & 4.86e-3 & 4.77e-3 & 4.77e-3 & 3.91e-3 & 3.58e-3 \\ 
\hline                                                                                               
5 & 3.8e-3 & 5.82e-3 & 5.82e-3 & 5.82e-3 & 5.82e-3 & 5.82e-3 & 5.82e-3 & 5.82e-3 \\ 
\hline                                                                                               
6 & 9.0e-4 & 6.65e-4 & 6.55e-4 & 5.90e-4 & 5.46e-4 & 5.51e-4 & 5.13e-4 & 4.70e-4 \\ 
\hline                                                                                               
7 & 5.9e-4 & 5.84e-4 & 5.84e-4 & 5.84e-4 & 7.17e-4 & 5.84e-4 & 5.84e-4 & 5.84e-4 \\ 
\hline                                                                                               
8 & 1.7e-4 & 7.34e-5 & 7.41e-5 & 6.36e-5 & 5.87e-5 & 7.78e-5 & 5.60e-5 & 4.38e-5 \\ 
\hline                                                                                               
9 & 1.2e-4 & 5.86e-5 & 5.86e-5 & 5.86e-5 & 5.86e-5 & 5.86e-5 & 5.86e-5 & 5.86e-5 \\ 
\hline                                                                                               
10 & 3.4e-5 & 8.22e-6 & 7.98e-6 & 6.75e-6 & 8.76e-6 & 6.68e-6 & 6.41e-6 & 5.36e-6 \\
\hline                                                                                               
\end{tabular}                                                                                        
\caption{Interpolation errors for $f_2$ (see Equation \ref{f1.eqn}).}
\label{table:interp2}  
\end{table}                
For the Runge-type function $f_2$, we see behavior more reflective of the Lebesgue constant of the resulting nodal set: the Face nodal distribution returns back a higher error than all other optimized node sets.  For $f_2$, the Interpolatory Warp and Blend nodes for the pyramid actually give back the lowest interpolation error for high $N$, despite their Lebesgue constant being larger than that of the Fekete nodes.  For both the Greedy and QR-based ``approximate Fekete'' points, the error appears to oscillate around the error of the true Fekete points.    

Since the Lebesgue constant is only an upper bound on the interpolation error, the reported numerical errors may behave better than the Lebesgue constants would indicate.  

\subsection{Optimization of the Duplex/Interpolatory Warp and Blend nodes}

The final step in construction of the original Warp and Blend nodes on the triangle was the addition of a quadratic variation to the interior blend, creating a one-parameter family of nodal distributions.  For the tetrahedron, a similar quadratic term is added, both to the edge blending functions and the face blending functions.  A 1D optimization problem is then solved for the value of this parameter which minimizes the Lebesgue constant of the resulting distribution.  

Both the Duplex and Interpolatory Warp and Blend procedures for the pyramid may also be optimized in a similar fashion.  For the Duplex construction, one such option is to optimize the edge/face blending parameter associated with the shared tetrahedral face inside the pyramid.  For the Interpolatory Warp and Blend procedure, we may also add a quadratic variation to the edge and face functions.\footnote{For edges associated with one or more triangular faces, we multiply the blending function by $1+\alpha v_c^2$, where for each bordering triangular face, $v_c$ is the vertex that does not lie on the edge.  For triangular faces, we multiply the blending function by $1+\alpha p_e^2$, where $p_e(r,s,t)$ is the plane which is zero on the triangular face and 1 at the opposite edge, while for the quadrilateral face, we multiply the blending function by $1 + \alpha v_5^2$, where $v_5$ is the vertex function at the tip of the pyramid.}  However, in our numerical experiments, neither optimization procedure improved the Lebesgue constant of the resulting nodal distribution significantly.  Table~\ref{table:origOptLeb} shows the difference in Lebesgue constant after optimization of $\alpha$ for the triangle, tetrahedron, Duplex pyramid, and Interpolatory Warp and Blend pyramid.\footnote{The unoptimized Duplex blending parameter is taken to be $\alpha$ instead of 0.} While the effect of optimization is noticable for both the 2D triangle to the 3D tetrahedron, the effect of optimization for both Warp and Blend pyramid nodes is much less pronounced.  

Since this optimization depends completely on the choice of quadratic variation, we are currently investigating other parametrizations of blending functions to improve optimization of Lebesgue constant.  


\begin{table}[!h]
\centering                                                          
\begin{tabular}{|c|c|c|c|c|c|c|c|c|}
\hline
& \multicolumn{2}{c|}{Triangle} &  \multicolumn{2}{c|}{Tetrahedron}  &  \multicolumn{2}{c|}{Duplex pyramid} &  \multicolumn{2}{c|}{Pyramid WB} \\ 
\hline
N  & No opt & Opt &  No opt & Opt  & No opt & Opt & No opt & Opt\\ 

\hhline{|=|=|=|=|=|=|=|=|=|}                                         
3 & 3.12 & 3.12 & 2.93 & 2.93 & 2.73 & 2.73 & 2.75 & 2.74 \\        
\hline                                                              
4 & 3.82 & 3.70 & 4.07 & 4.07 & 3.80 & 3.80 & 3.90 & 3.77 \\        
\hline                                                              
5 & 4.55 & 4.27 & 5.36 & 5.32 & 5.06 & 5.06 & 5.11 & 5.10 \\        
\hline                                                              
6 & 5.69 & 4.96 & 7.38 & 7.01 & 6.66 & 6.61 & 7.23 & 7.00 \\        
\hline                                                              
7 & 7.02 & 5.74 & 9.82 & 9.21 & 9.65 & 9.48 & 9.75 & 9.74 \\        
\hline                                                              
8 & 9.16 & 6.67 & 13.75 & 12.54 & 14.65 & 14.32 & 14.22 & 14.20 \\  
\hline                                                              
9 & 11.83 & 7.90 & 18.85 & 17.02 & 23.39 & 23.11 & 20.82 & 20.76 \\ 
\hline                                                              
10 & 16.06 & 9.36 & 27.02 & 24.40 & 40.12 & 39.73 & 32.16 & 32.00 \\
\hline                                                              
\end{tabular}                                                       
\caption{Original and optimized Lebesgue constants for the triangle, tetrahedron, Duplex pyramid, and Interpolatory Warp and Blend pyramid.  The left columns show the unoptimized Lebesgue constant, and the right columns show the Lebesgue constants after optimization of the construction in the previous column.}
\label{table:origOptLeb}
\end{table}              


\section{Conclusions}

We have compared several methods for the construction of nodal sets on pyramids for conforming finite element methods.  Both explicit and iterative procedures are considered --- previously existing algorithms for computing both Fekete and ``approximate Fekete'' points are adapted to the pyramid, and a new Duplex pyramid Warp and Blend procedure is introduced.  Furthermore, a new Interpolatory Warp and Blend procedure is developed and applied to the pyramid.  Similarly to the triangle, Warp and Blend-based nodal sets deliver lower Lebesgue constants for moderate values of $N$, while the iteratively determined Fekete nodes give lower Lebesgue constants for $N \geq 7$.  The condition number of the Vandermonde matrix and interpolation error for two functions are also given.  

A directory containing files for the relevant nodal distributions is available for download on Github at \url{https://github.com/tcew/nodes}.  Both \textsc{Matlab}\texttrademark \verb+.mat+ files and script files with which to generate the given nodal sets are provided.  





\bibliography{pyramids}{}
\bibliographystyle{plain}

\end{document}